\documentclass[11pt]{article}
\usepackage{graphicx}
\usepackage{amsmath,amssymb,amsthm,amsfonts}
\usepackage{amssymb}
\newtheorem{thm}{Theorem}[section]

\newtheorem{lem}{Lemma}[section]
\newtheorem{prop}[thm]{Proposition}

\theoremstyle{definition}

\usepackage{appendix}
\newtheorem{ex}{Example}[section]

\numberwithin{equation}{section}

\DeclareMathSymbol{\C}{\mathalpha}{AMSb}{"43}

\newcommand{\eps}{\varepsilon}

\newcommand{\lam}{\lambda}

\newcommand{\R}{{\mathbb{R}}}
\newcommand{\h}{{\mathcal{H}}}
\newcommand{\inte}{\int_{\mathbb{R}^2}}

\newcommand{\bsub}{\begin{subequations}}
\newcommand{\esub}{\end{subequations}$\!$}

\begin{document}

\title{Blow-up solutions for two coupled Gross-Pitaevskii equations
%for a  Nonlinear Schr\"{o}dinger System
with attractive interactions\thanks{Email: \texttt{yjguo@wipm.ac.cn}; \texttt{xyzeng@whut.edu.cn}; \texttt{hszhou@whut.edu.cn}}}
\author{Yujin Guo$^a$
,\ Xiaoyu Zeng$^b$\ \ and Huan-Song Zhou$^b$\thanks{Corresponding author.}\\
\small $^a$Wuhan Institute of Physics and Mathematics,
    Chinese Academy of Sciences,\\
    \small P.O. Box 71010, Wuhan 430071, P. R. China\\
 \small $^b$Department of Mathematics, Wuhan University of Technology, Wuhan 430070, China   }

%\author{Yujin Guo\thanks{Wuhan Institute of Physics and Mathematics,
%    Chinese Academy of Sciences, P.O. Box 71010, Wuhan 430071,
%    P. R. China.  Email: \texttt{yjguo@wipm.ac.cn}. Y. J. Guo is partially supported by NSFC grant No. 11322104, as well as the Scientific Research Foundation for the Returned Overseas Chinese Scholars, State Education Ministry.}
%,\, Xiaoyu Zeng\thanks{Wuhan Institute of Physics and Mathematics,
%    Chinese Academy of Sciences, P.O. Box 71010, Wuhan 430071,
%    P. R. China.  Email: \texttt{xyzeng@wipm.ac.cn}. } ,\,  and\, Huansong Zhou\thanks{Wuhan Institute of Physics and Mathematics,
%    Chinese Academy of Sciences, P.O. Box 71010, Wuhan 430071,
%    P. R. China.  Email: \texttt{hszhou@wipm.ac.cn}.  H. S. Zhou
%    is partially supported by NSFC grants No. 11171339 and 11271360.}
%}

\date{}

\smallbreak \maketitle

\begin{abstract}
The paper is concerned with a system of two coupled time-independent Gross-Pitaevskii equations in $\mathbb{R}^2$, which is used to model  two-component Bose-Einstein condensates with both attractive intraspecies and attractive interspecies interactions. This system is essentially an eigenvalue problem of a stationary nonlinear Schr\"odinger system in $\mathbb{R}^2$, solutions of the problem are obtained by seeking minimizers of the associated variational functional with constrained mass (i.e. $L^2-$norm constaints). Under certain type of trapping potentials $V_i(x)$ ($i=1,2$), the existence, non-existence and uniqueness of this kind of solutions are studied. Moreover, by establishing some delicate energy estimates, we show that each component of the solutions blows up at the same point (i.e., one of the global minima of $V_i(x)$) when the total interaction strength of intraspecies and interspecies goes to a critical value. An optimal blowing up rate  for the solutions of the system is also given.

\end{abstract}

\vskip 0.2truein



\vskip 0.2truein

%\tableofcontents

\section{Introduction}

In this paper, we study the following system of two coupled time-independent Gross-Pitaevskii equations:
\begin{equation}\label{sys}
\begin{cases}
 -\Delta u_{1}+V_1(x)u_{1}=\mu_{1} u_{1}+b_1
u_{1}^3 +\beta u_{2}^2u_{1}\quad \mbox{in}\,\ \R^2,\\
-\Delta u_{2}+V_2(x)u_{2}=\mu_{2} u_{2}+b_2 u_{2}^3 +\beta
u_{1}^2u_{2}\quad \mbox{in}\,\ \R^2.
\end{cases}
\end{equation}
 The system (\ref{sys}) arises in describing two-component Bose-Einstein condensates (BEC), where
$(V_1(x),V_2(x))$ is a certain type of trapping potentials, $(\mu_1,\mu_2)\in\R \times \R$ is the chemical potential, $b_i$ ($i=1,\,2$) and $\beta$ are the interaction strength of cold atoms inside each component (i.e. intraspecies) and  between two components (i.e. interspecies), respectively. Here $b_i>0$ (or $\beta>0$) represents the intraspecies (or interspecies) interaction is attractive, otherwise, it is repulsive. Much attention has been paid to the experimental studies of BECs since the BEC phenomena were successfully observed in 1995 in the pioneering experiments \cite{Anderson,D}. After that various BEC phenomena are observed in BEC experiments such as the symmetry
breaking, the collapse, the appearance of quantized vortices in
rotating traps, the phase segregation, etc., which inspired the theoretical investigation on BECs, specially on the Gross-Pitaevskii equations--the fundamental model of describing the BEC.
Multiple-component BECs can display
more interesting phenomena absent in single-component BEC.
In the past decade, under variant conditions on $V_i(x)$, $b_i$ and $\beta$, the analogues of (\ref{sys}) have been investigated widely, see e.g. \cite{Bao,BJS,CZ,CZ2,CZ3,CS,HME,IT,KL,LW2,MPS,Pom,R,Ti,ZZZ} and the references therein. In these mentioned papers, the authors are concerned with either the  semiclassical states of the system (i.e. replacing $-\Delta$ by $-\varepsilon\Delta$ in (\ref{sys})) \cite{LW2,MPS,IT,CZ2,Pom}, or the problem (\ref{sys}) with $V_i(x)$ being constants \cite{BJS,CZ,CZ3,LW,WW}, or the problem (\ref{sys}) in repulsive cases \cite{Bao,CS,HME,KL,R}, etc. The results in these papers are mainly on studying the existence of positive solutions for small $\varepsilon>0$, the location of the maximum point of the solutions and the behavior of the solutions with respect to $\varepsilon \rightarrow 0$ or the coupling parameter tends to infinity. To the authors' knowledge, there seems few results concerning the blowing-up analysis on the $L^2-$ normalized solutions of (\ref{sys}) with $b_i>0$.

In this paper, we are interested in dealing with the blowing-up properties of solutions for system (\ref{sys}) with both attractive intraspecies and interspecies interactions, that is, $b_i>0$ and $\beta>0$ (totally attractive case).  We leave the case of the attractive intraspecies (i.e. $b_i>0$) and the repulsive interspecies (i.e. $\beta <0$) to the companion work \cite{GZZ2}.
%The phenomena of BEC have been studied  intensively in physics  experiments, since its first realization in dilute gases  of cold  alkali atoms \cite{Anderson,D}.
%In many experiments, a large number of (bosonic)
%atoms are confined to a trap and cooled to very low temperatures,
%and the condensation of a large fraction of particles into the same
%one-particle state occurs under a certain critical temperature.
%These BECs display various interesting quantum phenomena \cite{Zwerger,Cooper,D,Fetter}, such as the symmetry breaking, the collapse, the appearance of quantized vortices in rotating traps, the effective lower dimensional behavior in strongly elongated traps, and etc. The system of  multiple-component BECs has recently been realized experimentally in \cite{HME,MMR,MBG,PPW} and references therein, and more interesting phenomena absent in  single-component BEC were observed. For example, two-component BECs were first realized in a magic trap in rubidium, after which two-component BECs have been realized in many kinds of configurations, cf. \cite{MBG}: a single isotope that is in two different hyperfine spin states, two different isotopes of the same atoms, or isotopes of two different atoms.
%The intraspecies interactions among the atoms of each component BEC can be either attractive or repulsive.
The case of repulsive intraspecies and interspecies interactions was studied recently in \cite{KL,R} and the references therein. When $b_i>0$ and $\beta>0$,
as mentioned in \cite{R}, (\ref{sys}) is  certainly a very different problem because the energy functional may not be bounded from below. In totally attractive case, we may expect from the single component BEC, see e.g., \cite{D,GS,GZZ}, that the  collapse still happens if the particle number increases beyond a critical value. In addition to the
intraspecies interaction among atoms in each component, there exist interspecies interactions
among the components for multiple-component BECs. Therefore, multiple-component BECs present more complicated characters than single component BEC, and the
corresponding analytic investigations are more challenging.

%Our main interest of the present paper  is to investigate two-component Bose-Einstein condensates with  attractive intraspecies and attractive interspecies interactions. In the companion work \cite{GZZ2}, we analyze two-component Bose-Einstein condensates with attractive intraspecies and  repulsive interspecies interactions.

It is well-known that the system (\ref{sys}) can be obtained from the associated time dependent nonlinear Schr\"{o}dinger equations if one seeks for  the following type standing-wave solutions
$$(\psi_1(x,t),\psi_2(x,t))=(u_1(x)e^{-\mathbf{i}\mu _1 t},u_2(x)e^{-\mathbf{i}\mu _2 t}),\,\ \mbox{where} \,\ \mathbf{i}^2=-1.$$
(\ref{sys}) is essentially an eigenvalue problem of a system of two stationary nonlinear Schr\"{o}dinger equations,
%Moreover, (\ref{sys}) is related to
which is also the system of Euler-Lagrange equations ($\mu_1, \mu_2$ are the Lagrange multipliers) of the following constrained minimization problem,
\begin{equation}\label{eq1.4}
\hat e(b_1,b_2,\beta):=\inf_{\{(u_1,u_2)\in \mathcal{M}\}} E_{b_1,b_2,\beta}(u_1,u_2),
\end{equation}
where $\mathcal{M}$  is  the so-called mass constraint
\begin{equation}\label{norm}
\mathcal{M}=\Big\{(u_1,u_2)\in\mathcal{X}:\, \inte|u_1|^2dx=\inte|u_2|^2dx=1\Big\}
\end{equation}
and $\mathcal{X}=\h_1\times \h_2$ with
\begin{equation}\begin{split}
   \h_i &= \Big \{u\in  H^1(\R ^2):\ \int _{\R ^2}  V_i(x)|u(x)|^2 dx<\infty \Big\},\\
   \|u\|_{_{\h_i}}&=\Big(\int _{\R ^2} \Big[|\nabla u|^2+ V_i(x)|u(x)|^2\Big] dx\Big)^{\frac{1}{2}}, \,\ \mbox{where}\,\ i=1,\,2.
   \end{split} \label{H}
\end{equation}
The energy functional $ E_{b_1,b_2,\beta}(u_1,u_2)$ is given by
\begin{equation}\label{f1}
\begin{split}
   E_{b_1,b_2,\beta}(u_1,u_2)=&\sum_{i=1}^2\int_{\R ^2} \Big(|\nabla
  u_i|^2+V_i(x)|u_i|^2-\frac{b_i}{2}|u_i|^4\Big)dx\\
  &-\beta\inte |u_1|^2|u_2|^2dx \,,\quad (u_1,u_2)\in\mathcal{X}.
\end{split}\end{equation}
%Essentially, by the analysis similar to that of \cite{GWZZ} where one-component minimization problem is considered,
In this paper, we only interested in the solutions of (\ref{sys}) with constrained mass, that is, the minimizers of (\ref{eq1.4}).
We mention that, different from the single component minimization problem (see, e.g. \cite{GWZZ}), in the two-component case it seems not clear whether a minimizer of (\ref{eq1.4}) is a ground state of (\ref{sys}). But we know  that if $(u_1,u_2)$ is a minimizer of (\ref{eq1.4}), then $(u_1,u_2)$ is a positive solution of (\ref{sys}) for some Lagrange multiplier $(\mu_1,\mu_2)\in \R\times\R$.
In this paper, $V_i(x)$ are  trapping potentials  of the following type
\begin{equation}\label{1.12}
V_i(x)\in L^\infty_{\rm
loc}(\R^2),\,\ \lim_{|x|\to\infty} V_i(x) = \infty\,\ \text{and}\,\ \inf_{x\in
\R^2} V_i(x) =0,\,\ i=1,\,2.
\end{equation}

Throughout the paper, we  assume that  both $\inf_{x\in \R^2} \big(V_1(x)+V_2(x)\big) $ and $\inf_{x\in \R^2} V_i(x)$ are attained. Since problem (\ref{eq1.4}) is invariant by adding suitable
constants to $V_i(x)$ and  we may simply assume that $\inf_{x\in \R^2} V_i(x) = 0$. Moreover, in what follows we assume that  $b_i>0$ ($i=1,\,2$) and $\beta >0$, and denote the norm of $L^p(\R^2)$ by $\|\cdot\|_p$ for $p\in(1,\infty)$.

We study problem (\ref{eq1.4}) which is motivated by the recent works \cite{Bao2,GS,GWZZ,GZZ}, etc., where the following single component  minimization problem was investigated:
\begin{equation}\label{1:two}
e_i(a):=\inf_{\{u\in\h_i, \|u\|_2^2=1\}} E_{a}^i(u),\,\ a>0,
\end{equation}
and the energy functional $E_{a}^i(u)$ satisfies
\begin{equation}\label{single}
E_{a}^i(u):=\int_{\R ^2} \big(|\nabla
  u(x)|^2+V_i(x)|u(x)|^2\big)dx-\frac{a}{2}\int_{\R ^2}|u(x)|^4dx, \,\  i=1 \,\  \mbox{or}\,\  2.
\end{equation}
Actually, it was proved in \cite{Bao2,GS} that (\ref{1:two}) admits minimizers if and only if $0<a<a^*:=\|Q\|^2_2$, where $Q=Q(|x|)>0$ denotes  %(cf.  \cite{GNN,K,Li,mcleod})
the unique positive solution of the following nonlinear scalar field
equation
\begin{equation}
-\Delta u+ u-u^3=0\  \mbox{  in } \  \R^2,\  \mbox{ where }\  u\in
H^1(\R ^2).  \label{Kwong}
\end{equation}
On the other hand, the
following Gagliardo-Nirenberg inequality   was shown in \cite{W} that
\begin{equation}\label{GNineq}
\inte |u(x)|^4 dx\le \frac 2 {\|Q\|_2^{2}} \inte |\nabla u(x) |^2dx
\inte |u(x)|^2dx ,\   \  u \in H^1(\R ^2),
\end{equation}
where the identity  is achieved at $u(x) = Q(|x|)$. One can note from
(\ref{Kwong}) that $Q(|x|)$  satisfies
\begin{equation}\label{1:id}
\inte |\nabla Q |^2dx  =\inte Q ^2dx=\frac{1}{2}\inte Q ^4dx ,
\end{equation}
see also Lemma 8.1.2 in \cite{C}. Moreover, we have
 \begin{equation}
Q(x) \, , \ |\nabla Q(x)| = O(|x|^{-\frac{1}{2}}e^{-|x|}) \,\
\text{as} \,\ |x|\to \infty,  \label{4:exp}
\end{equation}
see \cite[Proposition ~4.1]{GNN} for more details.

%The inequality (\ref{GNineq}) is used to discuss the existence and nonexistence of minimzers for the following single-component BEC
%\begin{equation}\label{1:two}
%e_i(a_i):=\inf_{\{u\in\h_i,\int_{\mathbb{R}^2} u^2dx=1\}}
%E_{a_i}(u),
%\end{equation}
%where the GP energy functional is of the form
%\begin{equation}
%E_{a_i}(u):=\int_{\R ^2} \big(|\nabla
%  u(x)|^2+V_i(x)|u(x)|^2\big)dx-\frac{a_i}{2}\int_{\R ^2}|u(x)|^4dx\, , \quad  u\in \h_i \,, \label{single}
%\end{equation}
%where $i=1,\,2$. More precisely, by applying the inequality (\ref{GNineq}) as
%well as the concentration analysis, the following existence and
%non-existence of minimizers are established in \cite{GS}:
%

 We now introduce briefly the main results of the present paper. Our first result is concerned with the following existence and nonexistence of minimizers for (\ref{eq1.4}).

\begin{thm}\label{thm1}
 Let $Q$ be the unique positive radial solution of
(\ref{Kwong}) and $a^*=\|Q\|_2^2$. If $V_i(x)$ satisfies (\ref{1.12}) for $i=1,\,2$. Then,
\begin{enumerate}
\item [\rm(i).] When $0<b_1<a^*$, $0<b_2<a^*$ and  $\beta<\sqrt {(a^*-b_1)(a^*-b_2)}$, problem (\ref{eq1.4}) has at least one minimizer.
%$\mathcal{O}(b_1,b_2,\beta)>1$

\item [ \rm(ii).] Either $b_1> a^*$ or $b_2>a^*$ or $\beta> \frac{a^*-b_1}{2}+\frac{a^*-b_2}{2}$, (\ref{eq1.4}) has no minimizer.
%$\mathcal{O}(b_1,b_2,\beta)<1$ and
\end{enumerate}
\end{thm}

In order to prove Theorem \ref{thm1}, in Section 2 we introduce the following minimization  problem
\begin{equation}\label{1:eq1.5}
\small\mathcal{O}(b_1,b_2,\beta):=\inf_{\tiny\big\{\begin{array}{c}
 u_i\in H^1(\R^2), \\
\|u_i\|_2^2=1,
\end{array}i=1,2\big\}
}\frac{\inte (|\nabla u_1|^2+|\nabla u_2|^2)dx}{\frac{b_1}{2}\inte |u_1|^4dx
+\frac{b_2}{2}\inte |u_2|^4dx+\beta\inte |u_1|^2|u_2|^2dx}.
\end{equation}
We shall prove in Proposition \ref{prop:2.1} that if $\mathcal{O}(b_1,b_2,\beta)>1$, then there exists at least one minimizer for (\ref{eq1.4}).  However, if $\mathcal{O}(b_1,b_2,\beta)<1$, then there is no minimizer for (\ref{eq1.4}). Therefore, by Proposition \ref{prop:2.1},  to establish Theorem \ref{thm1} it suffices to evaluate $\mathcal{O}(b_1,b_2,\beta)$, depending on the range of $(b_1,b_2,\beta)$. On the other hand, as shown in Lemma A.2 of the appendix, Theorem 1.1 can be proved alternatively by applying directly the Gagliardo-Nirenberg inequality (\ref{GNineq}) and some recaling techniques. We remark that  some results similar to  Theorem \ref{thm1} were also proved  in \cite[Theorem 2.6]{Bao} by this kind of ideas.

Theorem \ref{thm1} gives a complete classification of the  existence and nonexistence of minimizers for (\ref{eq1.4}), except that $(b_1,b_2,\beta)$ satisfies
\begin{equation}\label{eq1.13}
0<b_1 \leq a^*,\,\ 0<b_2\leq a^*\text{ and }\beta\in\big[\sqrt {(a^*-b_1)(a^*-b_2)}, \frac{a^*-b_1}{2}+\frac{a^*-b_2}{2}\big].
\end{equation}
When $(b_1,b_2,\beta)$ satisfies (\ref{eq1.13}), in general it is difficult to evaluate $\mathcal{O}(b_1,b_2,\beta)$ so that one cannot employ directly Proposition \ref{prop:2.1} to discuss the existence of minimizers for (\ref{eq1.4}). Therefore, some new ideas are needed to address this  case. Under some additional assumptions on $(b_1,b_2,\beta)$, in this paper we shall derive the following existence and nonexistence of minimizers for (\ref{eq1.4}).

\begin{thm}\label{thm2}
Under condition  (\ref{1.12}) and let $\underset{x\in \R^2}\inf\big(V_1(x)+V_2(x)\big)$ be attained.
   If $0<b_1\not=b_2<a^*$ satisfies additionally
  \begin{equation}\label{eq1.14}  |b_1-b_2|\leq 2\sqrt{(a^*-b_1)(a^*-b_2)},\end{equation} then there exists $\delta=\delta(b_1,b_2)\in(\sqrt{(a^*-b_1)(a^*-b_2)},\frac{2a^*-b_1-b_2}{2}]$ such that
$$\mbox{for any}\ \, \beta\in[\sqrt{(a^*-b_1)(a^*-b_2)}, \delta(b_1,b_2)\big),$$
there exists at least one minimizer for (\ref{eq1.4}).
\end{thm}

\begin{thm}\label{thm03}Suppose $V_i(x)$ ($i=1,2$) satisfies   (\ref{1.12}) and assume  $\underset{x\in \R^2}\inf\big(V_1(x)+V_2(x)\big)$ is attained.
  If $0<b_1=b_2<a^*$ and $\beta>0$ satisfy  (\ref{eq1.13}). Then
  \begin{enumerate}
    \item [\rm(i)] Problem (\ref{eq1.4}) has  no minimizer if $\underset{x\in \R^2}\inf\big(V_1(x)+V_2(x)\big)=0$.
    \item [\rm(ii)]Problem (\ref{eq1.4}) has  at least one minimizer if
    $$  \hat e(a^*-\beta,a^*-\beta,\beta)<\underset{x\in \R^2}\inf\big(V_1(x)+V_2(x)\big).$$
  \end{enumerate}

\end{thm}

We mention that, for $b_1=b_2\in (0,a^*)$, condition (\ref{eq1.13}) implies $\beta=a^*-b_1=a^*-b_2$, then the point $(b_1,b_2,\beta)$ must be located on the segment $(a^*-\beta,a^*-\beta,\beta)$ with $\beta\in (0,a^*)$.

In view of Theorems \ref{thm1} to \ref{thm03}, the existence and nonexistence of minimizers for (\ref{eq1.4}) still remain open for the case where $0<b_1\not=b_2 \leq a^*$ and $\beta >0$ is close to $\frac{2a^*-b_1-b_2}{2}$ from below. By applying Proposition \ref{prop:2.1}, we shall derive Theorem \ref{thm2}  through establishing the estimate $\mathcal{O}(b_1,b_2,\beta)>1$ under the additional assumption (\ref{eq1.14}).  As for the case where $0<b_1=b_2<a^*$, it follows from (\ref{eq7.2}) that $\mathcal{O}(a^*-\beta,a^*-\beta,\beta)=1$, and we shall give the proof of Theorem \ref{thm03} (ii) by applying Ekeland's variational principle \cite[Theorem  5.1]{S}.

The proof of Lemma \ref{lemA3} in the appendix implies that the following relationship  always holds
\begin{equation}\label{A:a*}
0\leq \hat  e(a^*-\beta,a^*-\beta,\beta)\leq \inf_{x\in\R^2}\big(V_1(x)+V_2(x)\big).
\end{equation}
In particular, the second inequality in (\ref{A:a*}) can be strict for certain potentials $V_1(x)$ and $V_2(x)$. Here is an example:

%The following example illustrates that the strict inequality of (\ref{A:a*}) can happen
% $$ \hat e(a^*-\beta,a^*-\beta,\beta)<\inf_{x\in \R^2}\big(V_1(x)+V_2(x)\big)$$
% for some $V_1(x)$ and $V_2(x)$:

\begin{ex}\label{E1.1}
For any given  points $x_1$ and $x_2$ in $\R^2$ satisfying $|x_1-x_2|>4$,   consider the function $0\leq\zeta_i(x)\in C^2_0(B_1(x_i))$ satisfying $\|\zeta_i\|_2=1$, where $i=1,\,2$, and define a positive constant $C_\zeta$ by
$$C_\zeta:=\sum_{i=1}^2\int_{\R ^2} \Big(|\nabla
  \zeta_i|^2-\frac{a^*-\beta}{2}|\zeta_i|^4\Big)dx-{\beta}\inte |\zeta_1|^2|\zeta_2|^2dx<\infty. $$
Let $0\leq V_i(x)\in C^2(\R^2; \R)$ satisfy $V_i(x)=0 $ in  $ B_1(x_i)$ and $V_i(x)\ge 2C_\zeta$ in  $ B_2^c(x_i)$ as well as $\lim _{|x|\to\infty}V_i(x)=\infty$, where $i=1,\,2$.
One can check that
$$ \hat e(a^*-\beta,a^*-\beta,\beta)\leq E_{a^*-\beta,a^*-\beta,\beta}(\zeta_1,\zeta_2)=C_\zeta<2C_\zeta\le \inf_{x\in\R^2} \big(V_1(x)+V_2(x)\big).$$
\end{ex}

Example \ref{E1.1} and Theorem \ref{thm03} (ii) show us  that  for any fixed $0<\beta <a^*$, there exists at least one minimizer of (\ref{eq1.4}) at $(b_1,b_2)=(a^*-\beta,a^*-\beta)$ for some suitable potentials $V_1(x)$ and $V_2(x)$. On the other hand, Theorem \ref{thm03} (i) gives the non-existence of minimizers for (\ref{eq1.4}) at $(b_1,b_2)=(a^*-\beta,a^*-\beta)$ and $0<\beta <a^*$ for the case where $\inf_{x\in \R^2}\big(V_1(x)+V_2(x)\big)=0$ is attained.
%\noindent
%{\bf Questions: when $  \hat e(a^*-\beta,a^*-\beta,\beta)=\inf_{x\in \R^2}\big(V_1(x)+V_2(x)\big)$, is there any minimizer for  (\ref{eq1.3}) at $(b_1,b_2)=(a^*-\beta,a^*-\beta)$? We guess NO. If the answer is NO, how about the limit behavior of minimizers as $(b_1,b_2)\nearrow (a^*-\beta,a^*-\beta)$?}\\

Without loss of generality, in the follows one may restrict the minimization of (\ref{eq1.4}) to non-negative vector functions $(u_1,u_2)$, since $E_{b_1, b_2,\beta}(u_1,u_2)\geq
E_{b_1, b_2,\beta}(|u_1|,|u_2|)$  for any $(u_1,u_2)\in \mathcal{X}$, due to the fact that $\nabla |u_i|\leq |\nabla u_i|$ a.e. in $\R^2$ ($i=1,\,2$). We next discuss the uniqueness of non-negative minimizers for (\ref{eq1.4}). Our recent results in \cite{GS,GZZ} show that the single component minimization problem (\ref{1:two}) has a unique non-negative minimizer when the parameter $a>0$ is suitably small, and however there may exist multiple non-negative minimizers for (\ref{1:two}) as  $a \nearrow a^*$  under some classes of trapping potentials.  We shall prove in Section 4 that a similar uniqueness result also holds for (\ref{eq1.4}) if $|(b_1,b_2,\beta)|$ is suitably small.

\begin{thm}\label{thm1.3}
If $V_i(x)$ satisfies (\ref{1.12}) for $i=1$ and $2$, then (\ref{eq1.4}) admits a unique non-negative minimizer if $|(b_1,b_2,\beta)|$  is suitably small.
\end{thm}

 A similar result on uniqueness for the single component  problem (\ref{1:two}) was proved in \cite{Schnee,M} by using the contracting map. However, this method seems not work for our problem. In this paper, we prove Theorem \ref{thm1.3} by employing  an implicit function theorem.

%\begin{rem} From our above theorem, we see that the convergence  rate of $e(a_1,a_2)$ depends only on the decay rate of the common zero points of $V_1(x)$ and $V_2(x)$ as $(a_1,a_2)\nearrow (a^*,a^*)$.  Nevertheless, for the single BEC,
%it was proved in \cite[Lemma 3]{GS} that  if
%\begin{equation}\label{1:pi}p_i:=\max\{p_{ij},j=1,\ldots,
%n_i\}\,,\quad i=1,\,2\,,\end{equation} then there exists $m,M>0$, independent of $a_i$, such that
%\begin{equation}\label{1:ei}
%m (a^*-a_i)^\frac{p_i}{p_i+2}\leq e_i(a_i)\leq M
%(a^*-a_i)^\frac{p_i}{p_i+2}\quad \text{for}\quad 0\leq a_i\leq
%a^*\,,\quad i=1,\,2\,.
%\end{equation}
%\end{rem}

For any fixed $0<\beta <a^*$, we finally focus on the limit behavior of the minimizers for (\ref{eq1.4}) as $(b_1,b_2)\nearrow (a^*-\beta,a^*-\beta)$, i.e., $(b_1+\beta,b_2+\beta)\nearrow (a^*,a^*)$ in the case that  there is no minimizer for (\ref{eq1.4}) at the threshold $(b_1,b_2)=(a^*-\beta,a^*-\beta)$. In view of Theorem \ref{thm03} (i), we shall consider the special case where $\underset{x\in \R^2}\inf\big(V_1(x)+V_2(x)\big)=0$, $i.e.,$ the minima of $V_1(x)$  coincide with those of $V_2(x)$. More precisely, we assume that for $V_i(x)$ ($i=1,\,2$) takes the form of
\begin{equation}\label{1:V}
\begin{split}
&V_i(x) =   h_i(x) \prod_{j=1}^{n_i} |x-x_{ij}|^{p_{ij}}\ \  \mbox{with}\ \  C < h_i(x) < 1/C \ \  \mbox{in} \ \R^2, \\
 &\mbox{and}\ \, h_i(x)\in C_{\rm loc}^\alpha(\R^2) \ \text{ for some } \ \alpha\in(0,1),
\end{split}\end{equation}
where $n_i\in \mathbb{N}$, $p_{ij}>0$,  $x_{ik}\not=x_{ij}$ for $k\not=j$, and  $\underset{x\to x_{ij}}\lim h_i(x)$ exists for
all $1\leq j\leq n_i$. Without loss of generality, we also
assume that there exists $1\leq l\leq \min\{n_1,n_2\}$ such that
\begin{equation}\label{1:V1}
\begin{split}
x_{1j}&=x_{2j},\quad \mbox{where}\quad j=1,\ldots,l;\\
 x_{1j_1}&\neq x_{2j_2}, \quad \text{where}\quad j_i\in\big\{l+1,\ldots n_i\big\} \ \, \mbox{and}\ \, i=1,\,2.
 \end{split}
\end{equation}
Note that (\ref{1:V1}) implies
\begin{equation}\label{1:com} \Lambda:=\big\{x\in \R^2:\,
V_1(x)=V_2(x)=0\big\}=\big\{x_{11}, x_{12},\cdots,x_{1l}\big\}.
\end{equation}
Define
\begin{equation}\label{1:def}
\bar p_j:=\min\big\{p_{1j},p_{2j}\big\},\ \, j=1,\ldots,l;\quad
p_0:=\max_{1\leq j\leq l}\min\big\{p_{1j}, p_{2j}\big\}=\max_{1\leq j\leq
l}\bar p_j ,
\end{equation}
so that
\begin{equation}\label{1:ba}\bar \Lambda:=\big\{x_{1j}:\,  \bar
p_j=p_0,j=1,\ldots,l\big\}\subset\Lambda .
\end{equation}
Let $ \gamma_j\in (0,\infty]$  be given by
\begin{equation}\label{def:li}
 \gamma_j := \lim_{x\to x_{1j}} \frac{V_1(x)+V_2(x)}{|x-x_{1j}|^{p_0}}, \quad 1\leq j\leq l.
\end{equation}
Note that $\gamma_j<\infty$ if and only if
$x_{1j}\in\bar \Lambda$, where $1\leq j\leq l $. Finally, define $\gamma = \min\big\{\gamma_1,\ldots
,\gamma_l\big\}$ so that the set
\begin{equation}\label{def:Z}
\mathcal{Z}:=\big\{x_{1j}:\  \gamma_j=\gamma,\ \, 1\leq j\leq
l\big\}\subset\bar\Lambda
\end{equation}
denotes the locations of the flattest global minima of
$V_1(x)+V_2(x)$. Using above notations,  our main results can be stated as follows.

\begin{thm}\label{thm1.4}
Assume that $0<\beta <a^*$ and  $V_i(x)$ satisfies (\ref{1:V})-(\ref{1:V1})  for $i=1, 2$. Let $(u_{b_1},u_{b_2})$ be a non-negative
minimizer of (\ref{eq1.4}) as $(b_1,b_2)\nearrow (a^*-\beta,a^*-\beta)$. Then, for any  sequence
$\{(b_{1k},b_{2k})\}$ satisfying $(b_{1k}, b_{2k})\nearrow (a^*-\beta, a^*-\beta)$ as
$k\to\infty$,   there exists a subsequence of $\{(b_{1k}, b_{2k})\}$, still denoted by
$\{(b_{1k}, b_{2k})\}$, such that, for $i=1,\,2$, each $u_{b_{ik}}$ has a unique global maximum point $x_{ik}\overset{k}\to\bar x_0$  for some $\bar x_0\in \mathcal{Z}$ and  \begin{equation}\label{11:rate1}
\lim _{k\to\infty}\frac{|x_{ik}-\bar x_0|}{\big(a^*-\frac{b_{1k}+b_{2k}+2\beta}{2}\big)^\frac{1}{p_0+2}}= 0.
\end{equation}  Moreover, for $i=1,\,2$,
\begin{equation*}
\lim _{k\to\infty} \Big(a^*-\frac{b_{1k}+b_{2k}+2\beta}{2}\Big)^\frac{1}{p_0+2} u_{b_{ik}}\Big(\big(a^*-\frac{b_{1k}+b_{2k}+2\beta}{2}\big)^\frac{1}{p_0+2}
x+x_{ik}\Big)=\frac{\lambda}{\|Q\|_{2}}Q(\lambda x)
 \end{equation*}
 strongly in $H^1(\R^2)$, where $\lam >0$ is given by
\begin{equation}\label{def1:lam}
  \lambda = \Big(
\frac{p_0\gamma}{4}\inte | x|^{p_0}Q^2(x)dx\Big)^\frac{1}{p_0+2}
\end{equation}
for $p_0>0$ and $\gamma >0$ defined in (\ref{1:def}) and (\ref{def:Z}), respectively.
\end{thm}

%\begin{rem}
%In Theorem \ref{4:thm3}, we only assume that $\beta>0$ fixed and $a_i=b_i+\beta\nearrow a^*$. This
%denotes that when $\beta>a^*$ is fixed and  the intraspecies are
%repulsive, i.e., $b_i\nearrow a^*-\beta<0$, our Theorem \ref{4:thm3} is still  applicable.
%\end{rem}

%The proof of Theorem \ref{4:thm3} is based on the following optimal  energy estimates
%\begin{equation}\label{AA:energy}
%\begin{split}
%&C_1\Big(a^*-\frac{b_{1}+b_{2}+2\beta}{2}\Big)^\frac{p_0}{p_0+2}\leq \hat e(b_1,b_2,\beta)\\
%&\leq C_2
%\Big(a^*-\frac{b_{1}+b_{2}+2\beta}{2}\Big)^\frac{p_0}{p_0+2}\quad
%\mbox{as}\quad
% (b_1,b_2)\nearrow (a^*-\beta, a^*-\beta)\,,
% \end{split}
%\end{equation}
%where the positive
%constants $C_1$ and $ C_2$ are independent of $b_1$ and $b_2$. Even though the optimal upper bound of (\ref{AA:energy}) can be proved in a similar way of \cite{GS}, the methods of \cite{GS} does not give the optimal lower bound of (\ref{AA:energy}). As analyzed in details in Section 3, we shall  introduce some new ideas and prove the optimal energy estimate (\ref{AA:energy}).

Theorems \ref{thm1.3} and  \ref{thm1.4} imply that the symmetry breaking occurs in
the minimizers of (\ref{eq1.4}). Actually, consider the trapping potentials $V_1$ and $V_2$ of the form
\[
V_1(x)=V_2(x) =    \prod_{j=1}^{l} |x-x_{j}|^{p}, \ \ p>0,
\]
where the points $x_{j}$ with $j=1,\cdots , l$ are arranged on the vortices  of a
regular polygon. Then there exist $0<a_{*}\leq a_{**}<a^*$ such that for $0<b_i+\beta<a_*$ ($i=1$ and $2$), the  functional
(\ref{eq1.4}) has a unique non-negative minimizer by Theorem \ref{thm1.3}, which has the same symmetry as that of $V_1(x)=V_2(x)$. However, when $a_{**}<b_i+\beta<a^*$ ($i=1,\,2$), we  obtain from  Theorem \ref{thm1.4} that (\ref{eq1.4}) possesses    (at least) $l$ different non-negative minimizers, and both
components of the minimizers concentrate at a zero point of $V_1(x)=V_2(x)$, which imply the symmetry breaking.  We note that the symmetry
breaking bifurcation of ground states for single nonlinear Schr\"odinger/Gross-Pitaevskii equations
has been studied in detail in the literature, see, e.g., \cite{J,K11,K08}.

%One can show that symmetry breaking can only occur for
%$a_i$ and $\beta$ sufficiently large. That is, for small enough
%$a_i$ and $\beta$, there is always a unique minimizer (up to
%multiplication by a constant phase), as in the case $a_i=\beta=0$.
%This can be proved by using the technique employed in the proof of
%Theorem 1.1 in \cite{M}.

%Second, the idea of proving the symmetry breaking, $i.e.$ concentration at a unique point, in \ref{GS} is only effective in dealing with finite global minima of $V(x)$.

This paper is organized as follows: in
Section 2 we first derive the crucial  Proposition \ref{prop:2.1} on the auxiliary minimization problem $\mathcal{O}(b_1,b_2,\beta)$, based on which we then complete the proof of Theorem \ref{thm1}. In Section 3 we focus on the proof of Theorems \ref{thm2} and \ref{thm03}. More exactly,   we first use Proposition \ref{prop:2.1} to prove Theorem \ref{thm2}, and   Theorem \ref{thm03} is then proved  by applying Ekeland's variational principle. Theorem \ref{thm1.3} is then proved in Section 4 to address the uniqueness of nonnegative minimizers
for (\ref{eq1.4}) as $|(b_1,b_2,\beta)|$ is suitably small. In Section 5 we shall  establish Proposition \ref{3:thm2} on optimal energy estimates of minimizers, upon which we finally complete in Section 6 the proof of Theorem \ref{thm1.4}. In the appendix, we give  an alternative proof of Theorem \ref{thm1} and prove a lemma as well which is used in Section 2.

\section{Existence   of Minimizers}

In this section, we  address the proof of Theorem \ref{thm1} on the existence of minimizers. We start with introducing the following  auxiliary minimization problem
\begin{equation}\label{2:eq1.5}
\small\mathcal{O}(b_1,b_2,\beta):=\inf_{\tiny\big\{\begin{array}{ll}
 u_i\in H^1(\R^2), \\
\|u_i\|^2_2=1,
\end{array}
i=1,\,2\big\}
}\frac{\inte (|\nabla u_1|^2+|\nabla u_2|^2)dx}{\frac{b_1}{2}\inte |u_1|^4dx
+\frac{b_2}{2}\inte |u_2|^4dx+\beta\inte |u_1|^2|u_2|^2dx}.
\end{equation}
By analyzing (\ref{2:eq1.5}), our first aim  is to derive the following proposition, which gives a criteria on the existence of minimizers for (\ref{eq1.4}) based on the value of  $\mathcal{O}(b_1,b_2,\beta)$.

\begin{prop}\label{prop:2.1}
Suppose that (\ref{1.12}) holds. Let $b_1$, $b_2$ and $\beta$ be positive. Then
\begin{enumerate}
\item [\rm (i)] (\ref{eq1.4}) has at least one minimizer if $\mathcal{O}(b_1,b_2,\beta)>1$.
\item [\rm(ii)]   (\ref{eq1.4}) has no minimizer if $\mathcal{O}(b_1,b_2,\beta)<1$.
\end{enumerate}
\end{prop}

%Our analysis of this section and the next one shows that if $\mathcal{O}(b_1,b_2,\beta)=1$, then the existence of minimizers for (\ref{eq1.4}) depends heavily on the comparison of global minima between $V_1(x)$ and $V_2(x)$.
To establish Proposition \ref{prop:2.1}, we need  the following compactness lemma.

\begin{lem}\label{2:lem1}
Suppose $V_i\in L_{loc}^\infty(\R^2)$ satisfies
$\lim_{|x|\to\infty}V_i(x)=\infty$, where $i=1,\,2$. Then the embedding $\mathcal{X}=\h_1\times \h_2\hookrightarrow
L^{q}(\R^2)\times L^{q}(\R^2)$ is compact for all $2\leq q<\infty$.
\end{lem}

 Since Lemma \ref{2:lem1} can be proved in a similar way to that of  \cite[Theorem  XIII.67]{RS} or \cite[Theorem 2.1]{BW}, we omit the  proof.\qed

%We next discuss the continuity of $\mathcal{O}(\cdot)$.

\begin{lem}\label{lemA2}
Let $\mathcal{O}(\cdot)$ be defined by (\ref{2:eq1.5}), then $\mathcal{O}(\cdot)$ is locally  Lipschitz  continuous in $\R^3_+$.
\end{lem}

\noindent\textbf{Proof.}  We first prove that
\begin{equation}\label{eq7.2}
\frac{a^*}{\max\{b_1+\beta,b_2+\beta\}}\leq \mathcal{O}(b_1,b_2,\beta) \leq \frac{2a^*}{b_1+b_2+2\beta}.
\end{equation}
Indeed, the upper bound of (\ref{eq7.2}) follows directly by taking $(\frac{Q(x)}{\|Q\|_2},\frac{Q(x)}{\|Q\|_2})$ as a trial function of (\ref{2:eq1.5}). On the other hand, for any $(u_1,u_2)\in H^1(\R^2)\times H^1(\R^2)$ satisfying $\inte |u_i|^2dx=1, i=1,\,2,$ it then follows from the
Gagliardo-Nirenberg inequality (\ref{GNineq}) that
\begin{align*}
&\frac{\inte (|\nabla u_1|^2+|\nabla u_2|^2)dx}{\frac{b_1}{2}\inte |u_1|^4dx
+\frac{b_2}{2}\inte |u_2|^4dx+\beta\inte |u_1|^2|u_2|^2dx}\\
\geq &\frac{\frac{a^*}{2}\inte(|u_1|^4+|u_1|^4)dx}{\frac{b_1+\beta}{2}
\inte |u_1|^4dx+\frac{b_2+\beta}{2}\inte |u_2|^4dx}\geq \frac{a^*}{\max\{b_1+\beta,b_2+\beta\}},
\end{align*}
which then gives the lower bound of (\ref{eq7.2}). Therefore,  (\ref{eq7.2}) is proved.

Consider $(b_1,b_2,\beta),(\tilde b_1, \tilde b_2,\tilde\beta)\in \R^3_+$,
 and let $\{(u_{1n},u_{2n})\}$ be a minimizing sequence of $\mathcal{O}(b_1,b_2,\beta)$.  Since  (\ref{2:eq1.5}) is invariant under the rescaling: $u(x)\mapsto \lambda u(\lambda x), \lambda>0$, one may assume that
\begin{equation}\label{eq7.3}
 \inte (|\nabla u_{1n}|^2+|\nabla u_{2n}|^2)dx=1\  \text{ for all } n\in\mathbb{N}^+.
\end{equation}
We then obtain that
$$\inte |u_{in}|^4dx\leq \frac{2}{a^*}\inte|\nabla u_{in}|^2\le \frac{2}{a^*},\, \ i=1,\,2, $$
and
\begin{align*}\inte |u_{1n}|^2|u_{2n}|^2dx&\leq \frac{1}{2}\inte (|u_{1n}|^4+|u_{2n}|^4)dx\\
&\leq
\frac{1}{a^*}\inte \big(|\nabla u_{1n}|^2+|\nabla u_{2n}|^2\big)dx=\frac{1}{a^*}.\end{align*}
Furthermore, we have
\begin{align*}
\frac{1}{\mathcal{O}(b_1,b_2,\beta)}
=&\lim_{n\to\infty}\Big[\frac{\frac{\tilde b_1}{2}\inte |u_{1n}|^4dx+\frac{\tilde b_2}{2}\inte |u_{2n}|^4dx+\tilde \beta\inte |u_{1n}|^2|u_{2n}|^2dx}{\inte (|\nabla u_{1n}|^2+|\nabla u_{2n}|^2)dx}\\
 &+\overset{2}{\sum_{i=1}}\frac{b_i-\tilde b_i}{2}\inte|u_{in}|^4dx+(\beta-\tilde \beta)\inte |u_{1n}|^2|u_{2n}|^2dx \Big]\\
\leq &\lim_{n\to\infty}\Big[\frac{1}{\mathcal{O}(\tilde b_1,\tilde b_2,\tilde\beta)}
+\overset{2}{\sum_{i=1}}\frac{|b_i-\tilde b_i|}{2}\inte|u_{in}|^4dx+|\beta-\tilde \beta|\inte |u_{1n}|^2|u_{2n}|^2dx\Big] \\
\leq &\frac{1}{\mathcal{O}(\tilde b_1,\tilde b_2,\tilde\beta)}
+\overset{2}{\sum_{i=1}}\frac{|b_i-\tilde b_i|}{a^*}+\frac{|\beta-\tilde \beta|}{a^*},
\end{align*}
i.e.,
\begin{equation}\label{eq7.4}
\frac{1}{\mathcal{O}(b_1,b_2,\beta)}-\frac{1}{\mathcal{O}(\tilde b_1,\tilde b_2,\tilde\beta)}\leq \frac{3}{a^*}\big|(b_1,b_2,\beta)-(\tilde b_1, \tilde b_2,\tilde\beta)\big|.
\end{equation}
Similarly, taking $\{(\tilde u_{1n},\tilde u_{2n})\}$ as a minimizing sequence of
$\mathcal{O}(\tilde b_1,\tilde b_2,\tilde\beta)$, and repeating the above argument, we know  that
$$\frac{1}{\mathcal{O}(\tilde b_1,\tilde b_2,\tilde\beta)}-\frac{1}{\mathcal{O}(b_1,b_2,\beta)}\leq \frac{3}{a^*}\big|(b_1,b_2,\beta)-(\tilde b_1, \tilde b_2,\tilde\beta)\big|,$$
The above estimates then yield that
\begin{align*}\Big|\frac{\mathcal{O}(b_1,b_2,\beta)-\mathcal{O}(\tilde b_1,\tilde b_2,\tilde\beta)}{\mathcal{O}(b_1,b_2,\beta)\mathcal{O}(\tilde b_1,\tilde b_2,\tilde\beta)}\Big|
&=\Big|\frac{1}{\mathcal{O}(b_1,b_2,\beta)}-\frac{1}{\mathcal{O}(\tilde b_1,\tilde b_2,\tilde\beta)}\Big|\\
&\leq \frac{3}{a^*}\big|(b_1,b_2,\beta)-(\tilde b_1, \tilde b_2,\tilde\beta)\big|.\end{align*}
By applying (\ref{eq7.2}), we therefore conclude that
$$\big|\mathcal{O}(b_1,b_2,\beta)-\mathcal{O}(\tilde b_1,\tilde b_2,\tilde\beta)\big|\leq \frac{12a^*}{(b_1+b_2+2\beta)(\tilde b_1+\tilde b_2+2\tilde\beta)}\big|(b_1,b_2,\beta)-(\tilde b_1, \tilde b_2,\tilde\beta)\big|,$$
which implies that  $\mathcal{O}(\cdot)$ is locally  Lipschitz  continuous in $\R^3_+$. \qed

With Lemma \ref{lemA2}, we now prove Proposition \ref{prop:2.1}.

\noindent\textbf{Proof of Proposition \ref{prop:2.1}.} i)  Let
$\{(u_{1n},u_{2n})\}\subset \mathcal{X}$ be a minimizing sequence
of problem (\ref{eq1.4}), $i.e.$,
$$\|u_{1n}\|_2^2=\|u_{2n}\|_2^2=1\quad \text{and}\quad
\lim_{n\to\infty}E_{b_1,b_2,\beta}(u_{1n},u_{2n})=\hat e(b_1,b_2,\beta).$$
It then follows from (\ref{2:eq1.5}) that
\begin{eqnarray}\label{eq2.1}
  E_{b_1,b_2,\beta}(u_{1n},u_{2n})\geq\sum_{i=1}^2\inte\Big[\big(1-\frac{1}{\mathcal{O}(b_1,b_2,\beta)}\big) |\nabla u_{in}|^2+V_i(x)|u_{in}|^2\Big] dx.
\end{eqnarray}
If $\mathcal{O}(b_1,b_2,\beta)>1$,   (\ref{eq2.1}) implies that $\{(u_{1n},u_{2n})\}\subset \mathcal{X}$ is  bounded  in  $n$. Thus, by the compactness of Lemma
\ref{2:lem1}, there exist a subsequence of $\{(u_{1n},u_{2n})\}$  and $(u_1,u_2)\in
\mathcal{X}$ such that
\begin{equation*}\begin{split}&(u_{1n},u_{2n})\overset{n}{\rightharpoonup}(u_1,u_2)\quad\text{weakly
in}\ \, \mathcal{X}\,,\\&
(u_{1n},u_{2n})\overset{n}{\to}(u_1,u_2)\quad\text{strongly in}\ \,
L^q(\R^2)\times L^q(\R^2)\ \text{for all }q\in[2,\infty). \end{split}\end{equation*}
Therefore,  $\|u_1\|_2^2=\|u_2\|_2^2=1$ and $
E_{b_1,b_2,\beta}(u_1,u_2)=\hat e(b_1,b_2,\beta)$. This proves part one of the proposition.

ii) Suppose now that $\mathcal{O}(b_1,b_2,\beta)<1$. One then may choose $(u_1,u_2)\in\mathcal{M}$ such that each $u_i\ (i=1,2 )$ has compact support in $\R^2$ and satisfies
\begin{equation}\label{eq2.2}
\frac{\inte (|\nabla u_1|^2+|\nabla u_2|^2)dx}{\frac{b_1}{2}\inte |u_1|^4dx+\frac{b_2}{2}\inte |u_2|^4dx+\beta\inte |u_1|^2|u_2|^2dx}\leq \delta:=\frac{1+\mathcal{O}(b_1,b_2,\beta)}{2}<1.
\end{equation}
 For $\lambda>0$, define \begin{equation}\label{eq2.3}\bar u_i(x)=\lambda u_i(\lambda x), \ \, i=1,\,2,\end{equation}
so that  $(\bar u_1, \bar u_2)\in\mathcal{M}$. Since $u_i(x)$ is compactly supported in $\R^2$ and $V_i(x)\in L_{\rm loc}^\infty(\R^2)$, there exists a positive constant $C$, independent of $\lambda>0$, such that for $\lambda\to\infty$,
\begin{equation}\label{eq2.4}
\inte V_i(x) |\bar u_i|^2dx=\inte V_i(\frac{x}{\lambda})|u_i|^2dx\leq C<\infty , \ \, i=1,\,2.
\end{equation}
On the other hand, by (\ref{eq2.2}) and (\ref{eq2.3}), we have
\begin{equation} \label{eq2.5}\begin{split}
  &\sum_{i=1}^2\int_{\R ^2} \Big(|\nabla
  \bar u_i|^2-\frac{b_i}{2}|\bar u_i|^4\Big)dx-\beta\inte |\bar u_1|^2|\bar u_2|^2dx \\
  =&\,\lambda^2\sum_{i=1}^2\int_{\R ^2} \Big(|\nabla
  u_i|^2-\frac{b_i}{2}|u_i|^4\Big)dx-\beta\lambda^2\inte |u_1|^2|u_2|^2dx \\
  \leq &\,\lambda^2(\delta-1)\Big(\frac{b_1}{2}\inte |u_1|^4dx+\frac{b_2}{2}\inte |u_2|^4dx+\beta\inte |u_1|^2|u_2|^2dx\Big) \\
\to & -\infty \,\ \text{ as }\,\ \lambda\to\infty.
\end{split}\end{equation}
It then follows from (\ref{eq2.4}) and (\ref{eq2.5}) that
$$\hat e(b_1,b_2,\beta)\leq E_{b_1,b_2,\beta}(\bar u_1, \bar u_2)\to-\infty \text{ as } \lambda\to\infty,$$
which implies that $\hat e(b_1,b_2,\beta)$ does not admit any minimizer. Proposition \ref{prop:2.1} is therefore established.
\qed

%\subsection{Proof of Theorem \ref{thm1}}
We end this section by proving  Theorem \ref{thm1}.

\noindent \noindent\textbf{Proof of Theorem \ref{thm1}.} \textbf{(i):} For any $(u_1,u_2)\in H^1(\R^2)\times H^1(\R^2) $ with $\inte |u_i|^2dx=1,\,i=1,\,2$, it follows from  (\ref{GNineq}) that
\begin{align*}
&\mathcal{O}(b_1,b_2,\beta)\\
\geq &\inf_{
\{\|u_i\|^2_2=1,i=1,2\}
}\frac{\frac{a^*}{2}\inte (| u_1|^4+|u_2|^4)dx}{\frac{b_1}{2}\inte |u_1|^4dx+\frac{b_2}{2}\inte |u_2|^4dx+\beta\big(\inte |u_1|^4dx\inte|u_2|^4dx\big)^\frac{1}{2}}\\
=&\inf_{
\{\|u_i\|_2^2=1,i=1,2\}
}\frac{\frac{a^*}{2}\Big(1
+\inte | u_2|^4dx\big/\inte|u_1|^4dx\Big)}{\frac{b_1}{2}+\frac{b_2}{2}\inte |u_2|^4dx\big/\inte|u_1|^4dx+\beta\big(\inte|u_2|^4dx\big/\inte|u_1|^4dx\big)^\frac{1}{2}}.
\end{align*}
Setting
\begin{equation}\label{eq2.6}
t:=\big(\inte|u_2|^4dx\big/\inte|u_1|^4dx\big)^\frac{1}{2}\in(0,\infty) \text{ and } f_{b_1,b_2,\beta}(t):=\frac{\frac{a^*}{2}(1+t^2)}{\frac{b_1}{2}+\frac{b_2}{2}t^2+\beta t},
\end{equation}
and then
 \begin{equation}\label{eq2.77}
 \mathcal{O}(b_1,b_2,\beta)\geq \inf_{t\in(0,\infty)} f_{b_1,b_2,\beta}(t).
\end{equation}
Since $0<b_i<a^*$ ($i=1,\,2$) and $\beta<\sqrt {(a^*-b_1)(a^*-b_2)}$,  standard calculations show that
$f_{b_1,b_2,\beta}(t)>1 \text{ for any }t\in(0,\infty)$, and also $$\underset{t\to0^+}\lim{f_{b_1,b_2,\beta}}(t)=\frac{a^*}{b_1}>1\ \, \text{ and }\,\  \underset{t\to\infty}\lim f_{b_1,b_2,\beta}(t)=\frac{a^*}{b_2}>1.$$
So, by the continuity  of $f_{b_1,b_2,\beta}(t)$ we obtain that $\underset{t\in(0,\infty)} \inf f_{b_1,b_2,\beta}(t)>1$. This estimate and (\ref{eq2.77}) then imply that
$\mathcal{O}(b_1,b_2,\beta)>1$, from which we conclude that (\ref{eq1.4}) has at least one minimizer by Proposition \ref{prop:2.1}(i).

\textbf{(ii):} Consider a   function  $0\le \varphi\in C_0^\infty(\R^2)$ satisfying $\inte |\varphi|^2dx=1$, and set \begin{equation}\label{eq2.6A}
u_\lambda(x)=\frac{\lambda}{\|Q\|_2}Q(\lambda x),\ \lambda>0,
\end{equation}
where $Q(x)$ is the unique radial positive solution  of  the   scalar field equation (\ref{Kwong}). It then follows from (\ref{1:id}) that
$$\inte |\nabla u_\lambda|^2dx=\frac{\lambda^2}{\|Q\|_2^2}\inte |\nabla Q|^2dx=\lambda^2,$$
as well as
$$ \inte |u_\lambda|^4dx=\frac{\lambda^2}{\|Q\|_2^4}\inte | Q|^4dx=\frac{2\lambda^2}{a^*}.$$

Suppose now that $b_1>a^*$. We then take $(u_\lambda, \varphi)$ as a trial function of $\mathcal{O}$ so that
\begin{align*}
\mathcal{O}(b_1,b_2,\beta)&\leq \frac{\inte \big(|\nabla u_\lambda|^2+|\nabla \varphi|^2\big)dx}{\frac{b_1}{2}\inte |u_\lambda|^4dx+\frac{b_2}{2}\inte |\varphi|^4dx+\beta\inte |u_\lambda|^2|\varphi|^2dx}\\
&=\frac{\lambda^2+\inte |\nabla \varphi|^2dx}{\frac{b_1}{a^*}\lambda^2+\frac{b_2}{2}\inte |\varphi|^4dx+\frac{\beta}{a^*}\inte |Q|^2|\varphi(\frac{x}{\lambda})|^2dx}\\
&\leq \frac{\lambda^2+\inte|\nabla \varphi|^2dx}{\frac{b_1}{a^*}\lambda^2}\to\frac{a^*}{b_1}<1\,\
\text{ as } \,\ \lambda\to\infty.\end{align*}
Thus, $\mathcal{O}(b_1,b_2,\beta)<1$ and it then follows from Proposition \ref{prop:2.1}(ii) that (\ref{eq1.4}) has no minimizer.
Similarly, if $b_2>a^*$, one can also obtain the nonexistence of minimizers for (\ref{eq1.4}).

Assume finally that $\beta> \frac{a^*-b_1}{2}+\frac{a^*-b_2}{2}$. In this case, take $(u_\lambda, u_\lambda)$ as a trial function  of $\mathcal{O}$, where $u_\lambda\ge 0$ is defined by (\ref{eq2.6A}). We then have
\begin{align*}
 \mathcal{O}(b_1,b_2,\beta)\leq \frac{2\inte |\nabla u_\lambda|^2dx}{\big(\frac{b_1}{2}+\frac{b_2}{2}+\beta\big)\inte |u_\lambda|^4dx}=\frac{a^*}{\frac{b_1}{2}+\frac{b_2}{2}+\beta}<1.
 \end{align*}
Hence, it follows again from Proposition \ref{prop:2.1}(ii) that (\ref{eq1.4}) has no minimizer.
This completes the proof of Theorem \ref{thm1}. \qed

\section{Further results on the existence of minimizers}
As discussed in the Introduction, our Theorem \ref{thm1} gives a complete classification of the  existence   of minimizers for (\ref{eq1.4}), except that $(b_1,b_2,\beta)$ satisfies
\[
0<b_1 \leq a^*,\,\ 0<b_2\leq a^*\text{ and }\beta\in\big[\sqrt {(a^*-b_1)(a^*-b_2)}, \frac{a^*-b_1}{2}+\frac{a^*-b_2}{2}\big].
\]
The aim of  this section is to prove Theorems \ref{thm2} and \ref{thm03}, which are concerned with the existence of minimizers for (\ref{eq1.4}) when $(b_1,b_2,\beta)$ lies in the above range. It turn out that such an existence depends on whether $0<b_1=b_2 \leq a^*$ or not.

Firstly, we shall make full use of  Proposition \ref{prop:2.1} to derive the existence of minimizers for (\ref{eq1.4}) in the case where $0<b_1\not=b_2<a^*$ and $\beta$ is close to $\sqrt{(a^*-b_1)(a^*-b_2)}$ from above. We start with the following lemma.

 \begin{lem}\label{33:lem3.1}
Let $0<b_1\not=b_2<a^*$ and $\beta=\sqrt{(a^*-b_1)(a^*-b_2)}$ such that
$$|b_1-b_2|\leq 2\beta.$$
If $\mathcal{O}\big(b_1,b_2,\sqrt{(a^*-b_1)(a^*-b_2)}\big)$ possesses a radially symmetric (about the origin) minimizing sequence $\{(u_{1n},u_{2n})\}\subset H_r^1(\R^2)\times H_r^1(\R^2)$   satisfying
{\small \begin{equation}\label{eq3.1}
0<C_1\leq \inte |\nabla u_{in}|^2dx\leq C_2<\infty, \text{ and }
0<C_1\leq \inte|u_{in}|^4dx\leq C_2<\infty,\,\ i=1, 2,
 \end{equation}}
and $C_1$, $C_2$ are independent of $n$, then we have
$$\mathcal{O}\big(b_1,b_2,\sqrt{(a^*-b_1)(a^*-b_2)}\big)>1.$$
\end{lem}

\noindent\textbf{Proof.} Since $\{u_{in}\subset H_r^1(\R^2)\}$ ($i=1,\,2$) is  radially symmetric, by (\ref{eq3.1}) and the compactness lemma of Strauss \cite{Str}, we deduce that there exists $u_i(x)\in H_r^1(\R^2)$ satisfying
\begin{equation}\label{eq3.02}\begin{split}
u_{in}&\overset{n}\rightharpoonup u_i \,\ \text{ weakly in }\,\  H^1(\R^2), \ \, i=1,\,2;\\
u_{in}&\overset{n}\to u_i\,\  \text{ strongly in }\,\ L^p(\R^2), \, \ \forall\ p\in(2,\infty),\ \, i=1,\,2.\end{split}
\end{equation}
Also, the assumption (\ref{eq3.1}) implies that
$$\inte |u_{i}|^4dx\geq C_1>0 \text{ and }u_i\not\equiv 0, \ i=1,\,2.$$

Since $0<b_1\not=b_2<a^*$, without loss of generality, we may assume that $b_1<b_2$. Then the assumption on $\beta$ can be simplified as
\begin{equation}\label{eq3.03}
0<b_1<b_2<a^* \ \text{ and } \ b_2\leq 2\beta +b_1, \, \text{where }\, \beta=\sqrt{(a^*-b_1)(a^*-b_2)}.
\end{equation}
Applying (\ref{GNineq}) and (\ref{eq3.02}), we have
\begin{equation}\label{eq3.4}\begin{split}
\mathcal{O}(b_1,b_2,\beta)&\geq \lim_{n\to\infty}\frac{\frac{a^*}{2}\inte (| u_{1n}|^4+|u_{2n}|^4)dx}{\frac{b_1}{2}\inte |u_{1n}|^4dx+\frac{b_2}{2}\inte |u_{2n}|^4dx+\beta\inte |u_{1n}|^2|u_{2n}|^2dx} \\
&=\frac{\frac{a^*}{2}\inte (| u_{1}|^4+|u_{2}|^4)dx}{\frac{b_1}{2}\inte |u_{1}|^4dx+\frac{b_2}{2}\inte |u_{2}|^4dx+\beta\inte |u_{1}|^2|u_{2}|^2dx} \\
&\geq
% \frac{\frac{a^*}{2}\inte (| u_{1}|^4+|u_{2}|^4)dx}{\frac{b_1}{2}\inte |u_{1}|^4dx+\frac{b_2}{2}\inte |u_{2}|^4dx+\beta\big(\inte |u_{1}|^4dx\cdot\inte|u_{2}|^4dx\big)^\frac{1}{2}}\\
% &=:
f_{b_1,b_2,\beta}(t_0), \,\ t_0:=\big(\inte|u_2|^4dx\big/\inte|u_1|^4dx\big)^\frac{1}{2}\in(0,\infty),
\end{split}
\end{equation}
where $f_{b_1,b_2,\beta}(t)$ is defined in (\ref{eq2.6}), and the equality of (\ref{eq3.4}) holds if and only if
\begin{equation}\label{eq3.6}
u_2^2(x)=\kappa u_1^2(x) \,\ \text{ for some }\, \ \kappa>0.
\end{equation}
Moreover, since $\beta=\sqrt{(a^*-b_1)(a^*-b_2)}$, it holds that
\begin{equation}\label{eq3.07}
f_{b_1,b_2,\beta}(t)\geq1,\,\  \forall\  t\in(0,\infty),
\end{equation}
and
\begin{equation}\label{eq3.8}
f_{b_1,b_2,\beta}(t)=1 \,\ \Leftrightarrow\,\  t=t_1:=\sqrt\frac{a^*-b_1}{a^*-b_2}.
\end{equation}
We thus deduce from (\ref{eq3.4}) and (\ref{eq3.07}) that $\mathcal{O}\big(b_1,b_2,\sqrt{(a^*-b_1)(a^*-b_2)}\big)\geq 1$.

We claim that   $\mathcal{O}\big(b_1,b_2,\sqrt{(a^*-b_1)(a^*-b_2)}\big)>1$. Otherwise, if
\begin{equation}\label{eq3.09}\mathcal{O}\big(b_1,b_2,\sqrt{(a^*-b_1)(a^*-b_2)}\big)=1,
\end{equation}
then (\ref{eq3.6})  holds for $t_0=t_1$, where $t_0$ and $t_1$ are given by (\ref{eq3.4}) and (\ref{eq3.8}), respectively.  Thus, we  have
\begin{equation}\label{eq3.9}
u_2^2(x)=\kappa u_1^2(x),\,\ \text{where}\,\ \kappa=t_0=\sqrt\frac{a^*-b_1}{a^*-b_2}>1.
\end{equation}
Together with (\ref{eq3.02}), this implies  that
\begin{align}
\mathcal{O}(b_1,b_2,\beta)&\geq \frac{\inte (| \nabla u_{1}|^2+|\nabla u_{2}|^2)dx}{\frac{b_1}{2}\inte |u_{1}|^4dx+\frac{b_2}{2}\inte |u_{2}|^4dx+\beta\inte |u_{1}|^2|u_{2}|^2dx}\nonumber\\
&=\frac{1+\kappa}{\frac{b_1}{2}+\frac{b_2}{2}\kappa^2+\beta\kappa}\cdot\frac{\inte |\nabla u_{1}|^2dx}{\inte |u_{1}|^4dx}.\label{eq3.10}
\end{align}

On the other hand, define
\begin{equation}\label{eq3.11}
\tilde u_i(x)=\frac{1}{\sqrt\lambda_i}u_i(x)\,\ \text{where}\,\  \lambda_i:=\inte |u_i|^2dx\leq \lim_{n\to\infty}\inte |u_{in}|^2dx=1,
\end{equation}
so that $\inte|\tilde u_i|^2dx=1$, where $i=1,\,2$. Note also
 from (\ref{eq3.9}) that
\begin{equation}\label{eq3.12}
\lambda_2=\kappa\lambda_1\leq1.
\end{equation}
Therefore, by the definition of $\mathcal{O}(\cdot)$, we deduce from (\ref{eq3.9}), (\ref{eq3.11}) and (\ref{eq3.12}) that
\begin{align}
\mathcal{O}(b_1,b_2,\beta)&\leq \frac{\inte (|\nabla\tilde u_{1}|^2+|\nabla \tilde u_{2}|^2)dx}{\frac{b_1}{2}\inte |\tilde u_{1}|^4dx+\frac{b_2}{2}\inte |\tilde u_{2}|^4dx+\beta\inte |\tilde u_{1}|^2|\tilde u_{2}|^2dx}\nonumber\\
&=\frac{2\lambda_1}{\frac{b_1}{2}+\frac{b_2}{2}+\beta}\cdot \frac{\inte | \nabla u_{1}|^2dx}{\inte |u_{1}|^4dx}\nonumber
\leq \frac{\frac2\kappa}{\frac{b_1}{2}+\frac{b_2}{2}+\beta}\cdot\frac{\inte |\nabla u_{1}|^2dx}{\inte |u_{1}|^4dx}.\nonumber
\end{align}
It then follows from  (\ref{eq3.10})  that
$$\frac{1+\kappa}{\frac{b_1}{2}+\frac{b_2}{2}\kappa^2+\beta\kappa}\leq \frac{\frac2\kappa}{\frac{b_1}{2}+\frac{b_2}{2}+\beta},$$
i.e.,
$$\frac{\kappa+\kappa^2}{\frac{b_1}{2}+\frac{b_2}{2}\kappa^2+\beta\kappa}\leq \frac{2}{\frac{b_1}{2}+\frac{b_2}{2}+\beta}, \ \text{where }\ \kappa=\sqrt\frac{a^*-b_1}{a^*-b_2}>1,$$
which however  contradicts  Lemma \ref{lemA} in the Appendix. Thus (\ref{eq3.09}) cannot occur, then $\mathcal{O}\big(b_1,b_2,\sqrt{(a^*-b_1)(a^*-b_2)}\big)>1$. \qed

With Lemma \ref{33:lem3.1} and Proposition \ref{prop:2.1}, we can prove now  Theorem \ref{thm2}.

\noindent \noindent\textbf{Proof of Theorem \ref{thm2}.}
Let $\{(u_{1n},u_{2n})\}$ be a minimizing sequence of $\mathcal{O}(b_1,b_2,\beta)$. By the Schwarz symmetrization of $\{(u_{1n},u_{2n})\}$, one may assume that
\begin{equation}\label{eq2.211}
u_{in}(x)=u_{in}(|x|)\geq0,\, \ i=1,\,2.\end{equation}
Moreover, since  the problem (\ref{2:eq1.5}) is invariant under the rescaling: $u(x)\mapsto \lambda u(\lambda x)$ for $\lambda>0$, one can also assume that
\begin{equation}\label{eq3.2}
 \inte (|\nabla u_{1n}|^2+|\nabla u_{2n}|^2)dx=1\  \text{ for all } n\in\mathbb{N}^+.
\end{equation}
By the Gagliardo-Nirenberg inequality (\ref{GNineq}), we then obtain that $\inte |u_{in}|^4dx$ ($i=1,\,2$) is bounded uniformly, $i.e.,$
\begin{equation}\label{eq3.3}
0<C_1\leq \frac{b_1}{2}\inte |u_{1n}|^4dx+\frac{b_2}{2}\inte |u_{2n}|^4dx+\beta\inte |u_{1n}|^2|u_{2n}|^2dx\leq C_2<\infty.
\end{equation}
Under the assumption (\ref{eq1.14}), we  claim that
\begin{equation}\label{eq3.17}
\mathcal{O}\big(b_1,b_2,\sqrt{(a^*-b_1)(a^*-b_2)}\big)>1.
\end{equation}
We  prove (\ref{eq3.17}) by considering separately the following two cases.
 \vskip 0.1truein
 \noindent \textbf{Case 1.} If  \begin{equation}
 \inte|u_{in}|^4dx\to0\ \ \text{ as }\,\ n\to\infty, \ \, \text{ where }\, \ i=1 \text{ or }2.
\end{equation}
Without loss of generality, we  assume that
 $$\inte|u_{1n}|^4dx\to0\ \ \text{ as }\,\ n\to\infty.$$
It then follows from (\ref{eq3.3}) that
$$\inte |u_{1n}|^2|u_{2n}|^2dx\overset{n}\to0\,\ \text{ and }\,\ \inte |u_{2n}|^4dx\geq C>0.$$
By the assumption $0<b_2<a^*$, we then have
\begin{align}\mathcal{O}(b_1,b_2,\beta)
&=\lim_{n\to\infty} \frac{\inte \big(|\nabla u_{1n}|^2+|\nabla u_{2n}|^2\big)dx}{\frac{b_2}{2}\inte |u_{2n}|^4dx+o(1)}\nonumber\\
&\geq \lim_{n\to\infty} \frac{\inte |\nabla u_{2n}|^2dx}{\frac{b_2}{2}\inte |u_{2n}|^4dx+o(1)}\geq \frac{a^*}{b_2}>1,\label{eq3.18}
\end{align}
where  (\ref{GNineq}) is used.
Thus, (\ref{eq3.17}) follows immediately from (\ref{eq3.18}) with $\beta=\sqrt{(a^*-b_1)(a^*-b_2)}$.
 \vskip 0.1truein
\noindent \textbf{Case 2.} If
\begin{equation}\label{eq3.19}
\inte|u_{in}|^4dx\geq C>0 \ \ \text{for }\ i=1\,\text{ and }\, 2.
\end{equation}
In this case, applying (\ref{GNineq}) and (\ref{eq3.2}), we deduce that there exist positive constants $C_3$ and $C_4$, independent of $n$, such that
\begin{equation*}\begin{split}
0<&C_3\leq \inte |\nabla u_{in}|^2dx\leq C_4<\infty,\ \, i=1,\, 2,\\
0<&C_3\leq \inte|u_{in}|^4dx\leq C_4<\infty,\ \quad i=1,\, 2.
\end{split}\end{equation*}
Under the assumption (\ref{eq1.14}), we then conclude from (\ref{eq2.211}) and  Lemma \ref{33:lem3.1} that the estimate (\ref{eq3.17}) also holds. So, the claim is proved.

Combining (\ref{eq3.17}) and Lemma \ref{lemA2},  we then derive that there exists a constant $\delta $ satisfying $\delta>\sqrt{(a^*-b_1)(a^*-b_2)}$  such that
$$\mathcal{O}(b_1,b_2,\beta)>1 \ \,\text{ for all }\,\ \beta\in[\sqrt{(a^*-b_1)(a^*-b_2)}, \delta\big).$$
By Proposition \ref{prop:2.1} we therefore conclude that (\ref{eq1.4}) has at least one minimizer.
\qed\\

%\subsection{Proof of Theorem \ref{thm2} (II)}

Next,  we turn to proving Theorem \ref{thm03}, that is, the existence  and non-existence of minimizers for (\ref{eq1.4}) in the case where $(b_1,b_2,\beta)$ satisfies $0<b_1=b_2<a^*$ and (\ref{eq1.13}). These assumptions on $(b_1,b_2,\beta)$ imply  that $(b_1,b_2,\beta)$ lies on the segment:  $(b_1,b_2,\beta)=(a^*-\beta,a^*-\beta,\beta)$ with $0<\beta <a^*$.  %In order to simplify the notations, as stated in the introduction, we equivalently consider the minimization problem (\ref{eq1.3}) at the  threshold $(a_1,a_2)=(a^*,a^*)$.
In this case, for simplicity,  we rewrite the functional (\ref{f1})  as
\begin{equation}\label{f2}
\begin{split}
  E_{a^*,a^*}(u_1,u_2):=&\sum_{i=1}^2\int_{\R ^2} \Big(|\nabla
  u_i|^2+V_i(x)|u_i|^2-\frac{a^*}{2}|u_i|^4\Big)dx\\
  &+\frac{\beta}{2}\inte \big(|u_1|^2-|u_2|^2\big)^2dx \,,\quad\mbox{where}\quad  (u_1,u_2)\in\mathcal{X}\,.
\end{split}\end{equation}
Note from (\ref{2:limt2}) and (\ref{2:limt6}) in the Appendix  that
$$0\leq \hat e(a^*-\beta,a^*-\beta,\beta)\leq \inf_{x\in\R^2}\big(V_1(x)+V_2(x)\big).$$
If $\inf_{x\in\R^2}\big(V_1(x)+V_2(x)\big)=0$ is attained,  Theorem \ref{thm03}(i) shows  that (\ref{eq1.4}) does not admit any minimizer.  Surprisingly,  Theorem \ref{thm03}(ii) however shows that there exists at least one minimizer for (\ref{eq1.4})
%(\ref{eq1.3})
if
 $$0\le \hat e(a^*-\beta,a^*-\beta,\beta)< \inf_{x\in\R^2}\big(V_1(x)+V_2(x)\big).$$
 As illustrated in Example 1.1, the above condition does hold for some potentials  $V_1(x)$ and $V_2(x)$.

 Finally, we give the proof of Theorem \ref{thm03}.

%
%By considering the functional $E_{a^*,a^*}(u_1,u_2)$ of the form (\ref{f2}), in the following we address the proof of Theorem \ref{thm2} (II).\\

\noindent{\textbf{Proof of Theorem \ref{thm03} (i):}}  In order to prove part (i), we assume  that   $\inf_{x\in \R^2}\big(V_1(x)+V_2(x)\big)=0$, i.e.,  there exists  $x_0\in \R^2$ such that  $V_1(x_0)=V_2(x_0)=0$. Since  $b_1=b_2=a^*-\beta$, we take $\phi >0$ as in (\ref{2:trial}) and use $(\phi, \phi)$ with $\bar x_0 =x_0$ as a trial  function of $\hat e(a^*-\beta,a^*-\beta,\beta)$. It then follows from (\ref{2:limt2}) and (\ref{2:limt6}) that
\begin{equation}\label{5:ea*}
0\le \hat e(a^*-\beta,a^*-\beta,\beta)\leq \lim_{\tau\to \infty}E_{a^*,a^*}(\phi,\phi)=V_1( x_0)+V_2( x_0)=0,
\end{equation}
$i.e.,$  $\hat e(a^*-\beta,a^*-\beta,\beta)=0$. Suppose now  there exists a
minimizer $(\hat u_1,\hat u_2)\in\mathcal{M}$ for  $\hat e(a^*-\beta,a^*-\beta,\beta)$. As pointed out in the Introduction, we can assume $(\hat u_1,\hat u_2)$ to be non-negative. It then follows from (\ref{f2}) that $\hat u_1\equiv \hat u_2\ge 0$ in $\R^2$, and
\begin{equation*}
\int_{\R^2}|\nabla \hat u_1|^2dx=\frac{1}{2}\int_{\R^2}|\hat u_1|^4dx\,\ {\rm
and}\,\ \int_{\R^2}V_1(x)|\hat u_1|^2dx=0.
\end{equation*}
This is a contradiction, since the first equality implies that $\hat u_1(x)$
is equal to (up to translation) $Q(x)$, but the second equality yields
that $\hat u_1(x)$ has compact support. Therefore, the conclusion (i) of Theorem \ref{thm03}  is proved.

%Moreover, if $V_1(x_0)=V_2(x_0)=0$ for some $x_0\in \R^2$, then
%$\lim_{(a_1,a_2)\nearrow (a^*,a^*)} e(a_1, a_2) = e(a^*, a^*) = 0$.
%This follows easily from (\ref{2:limt2}) and
%(\ref{2:limt6}), by first taking $(a_1,a_2)\nearrow(a^*,a^*)$ and
%then letting $\tau\to \infty$. This implies that
%$\limsup_{(a_1,a_2)\nearrow(a^*,a^*)} e(a_1,a_2)\leq V_1(\bar
%x_0)+V_2(\bar x_0)$ which yields the result after taking the infimum over $\bar x_0$.

\textbf{(ii):} In this  case, we have  $$0\le \hat e(a^*-\beta,a^*-\beta,\beta)<
\inf_{x\in\R^2}\big(V_1(x)+V_2(x)\big).$$
For $\mathcal{M}$ defined by (\ref{norm}), we introduce
$$d(\vec u,\vec v):=\|\vec u-\vec v\|_\mathcal{X}, \quad   \vec u,\, \vec v\in \mathcal{M}\,,$$ where
$$\|\vec u\|_\mathcal{X}=\big(\|u_1\|_{\mathcal{H}_1}^2+\|u_2\|_{\mathcal{H}_2}^2\big)^\frac{1}{2},\quad \vec u=(u_1,u_2)\in \mathcal{X}\,.$$
It is easy to check that $(\mathcal{M},d)$ is a complete distance space. Hence, by Ekeland's variational principle \cite[Theorem  5.1]{S}, there exists a minimizing sequence  $\{\vec u_n=(u_{1n},u_{2n})\}\subset \mathcal{M}$ of $\hat e(a^*-\beta,a^*-\beta,\beta)$ such that
\begin{eqnarray}\label{eqE1}
&&\hat e(a^*-\beta,a^*-\beta,\beta)\leq E_{a^*,a^*}(\vec u_n)\leq \hat e(a^*-\beta,a^*-\beta,\beta)+\frac{1}{n}\,,\\
&&E_{a^*,a^*}(\vec v)\geq E_{a^*,a^*}(\vec u_n)-\frac{1}{n}\|\vec u_n-\vec v\|_\mathcal{X} \quad \text{for}\quad \vec v \in \mathcal{M}.\label{eqE2}
\end{eqnarray}
%As discussed in the Introduction, without loss of generality we may assume that the minimizing sequence is non-negative.

Due to the compactness of Lemma \ref{2:lem1}, in order to show that there exists a minimizer for $\hat e(a^*-\beta,a^*-\beta,\beta)$, it suffices to prove that  $\{\vec u_n=(u_{1n},u_{2n})\}$ is bounded in $\mathcal{X}$ uniformly w.r.t. $n$. We argue by contradiction. If $\{\vec u_n=(u_{1n},u_{2n})\}$ is unbounded in $\mathcal{X}$,  then there exists a subsequence of $\{\vec u_n\}$, still denoted by  $\{\vec u_n\}$, such that $\|\vec u_n\|_\mathcal{X}\xrightarrow{n}\infty$. By
Gagliardo-Nirenberg inequality, we deduce from  (\ref{eqE1}) that
\begin{equation}\label{5:vi}
\sum_{i=1}^2\inte V_i(x)|u_{in}|^2dx\leq E_{a^*,a^*}(\vec u_n)\leq \hat e(a^*-\beta,a^*-\beta,\beta)+\frac{1}{n}.
\end{equation}
Hence,
\begin{equation}\label{eqE4}\inte|\nabla u_{1n}|^2+|\nabla u_{2n}|^2dx\xrightarrow{n}\infty .\end{equation}
We now claim that
\begin{eqnarray}\label{5:com1}
&&\inte|\nabla u_{in}|^2dx\sim \frac{a^*}{2}\inte| u_{in}|^4dx\xrightarrow{n}\infty,\quad i=1,\, 2,\\
&&\inte| u_{1n}|^4dx\bigg/\inte| u_{2n}|^4dx\xrightarrow{n}1.\label{5:com2}
\end{eqnarray}
Indeed, by (\ref{eqE4}),  we  may assume that $\inte|\nabla u_{1n}|^2dx\xrightarrow{n}\infty$. Note from (\ref{eqE1}) that
\begin{equation}\label{5:ui}0\leq \inte|\nabla u_{in}|^2dx-\frac{a^*}{2}\inte| u_{in}|^4dx\leq \hat e(a^*-\beta,a^*-\beta,\beta)+\frac{1}{n},\ \,\mbox{for}\,\ i=1,\, 2.
\end{equation}
This implies that  \begin{equation}\label{5:u1}
\frac{a^*}{2}\inte| u_{1n}|^4dx\xrightarrow{n}\infty\quad \text{and}\,\ \inte|\nabla u_{1n}|^2dx\bigg/\frac{a^*}{2}\inte| u_{1n}|^4dx\xrightarrow{n}1. \end{equation}
On the other hand, (\ref{eqE1}) also yields that
\begin{equation}\label{5:eq}
\beta\inte\big(|u_{1n}|^2-|u_{2n}|^2\big)^2dx\leq \hat e(a^*-\beta,a^*-\beta,\beta)+\frac{1}{n}.
\end{equation}
Then,
\begin{align*}\begin{split}
&\bigg|\int_{\R^2}|u_{a_1}|^4dx-
\int_{\R^2}|u_{a_2}|^4dx\bigg|
\leq
\Big(\int_{\R^2}\big(|u_{a_1}|^2-|u_{a_2}|^2\big)^2dx\Big)^
\frac{1}{2}\Big(\int_{\R^2}\big(|u_{a_1}|^2+|u_{a_2}|^2\big)^2dx\Big)^\frac{1}{2}\\
&\leq
\left\{\frac{1}{\beta}\big(\hat e(a^*-\beta,a^*-\beta,\beta)+\frac{1}{n}\big)\right\}^{\frac{1}{2}}
\bigg[\Big(\int_{\R^2}|u_{a_1}|^4dx\Big)^\frac{1}{2}+\Big(\int_{\R^2}|u_{a_2}|
^4dx\Big)^\frac{1}{2}\bigg]
\end{split}\end{align*}
Together with (\ref{5:u1}), we thus derive that
\begin{equation}\label{5:u12}
\inte| u_{2n}|^4dx\xrightarrow{n}\infty \quad \text{and}\,\ \inte| u_{1n}|^4dx\bigg/\inte| u_{2n}|^4dx\xrightarrow{n}1.
\end{equation}
This estimate and (\ref{5:ui}) then imply that
\begin{equation}\label{5:u2}\inte| \nabla u_{2n}|^2dx\xrightarrow{n}\infty\quad \text{and}\,\ \inte|\nabla u_{2n}|^2dx\bigg/\frac{a^*}{2}\inte| u_{2n}|^4dx\xrightarrow{n}1. \end{equation}
Therefore,  (\ref{5:com1}) and (\ref{5:com2})  follow from   (\ref{5:u1}), (\ref{5:u12}) and (\ref{5:u2}).

Define now $$\epsilon_n^{-2}:=\inte |u_{1n}|^4dx.$$
Similar to Lemma \ref{le3.3} (i) in  Section \ref{Sec5}, there exists a sequence $\{y_{\epsilon_n}\}\subset\R^2$ as well as positive constants $R_0$ and $\eta$ such that
\begin{equation}\label{5:scal}w_{in}(x):=\epsilon_n u_{in}(\epsilon_n x+\epsilon_n y_{\epsilon_n}),\,\ i=1,\, 2,\end{equation} satisfies
\begin{equation}\label{5:R}
\int_{B_{R_0}(0)}|w_{in}|^2dx>\eta>0,\,\ i=1,\, 2.
\end{equation}
Recall from (\ref{5:eq}) that
\begin{equation*}
\inte\big(|w_{1n}|^2-|w_{2n}|^2\big)^2dx=\epsilon_n^2 \inte\big(|u_{1n}|^2-|u_{2n}|^2\big)^2dx\xrightarrow{n}0,
\end{equation*}
which implies that
\begin{equation}\label{5:L2}
w_{1n}^2-w_{2n}^2\xrightarrow{n}0 \,\ \text{in}\ L^2(\R^2)\,\ \text{and}\,\
w_{1n}^2-w_{2n}^2\xrightarrow{n}0 \,\ \text{a.e. in }\, \R^2.
\end{equation}
Moreover, since  $\lim_{|x|\to\infty}V_i(x)= \infty$  ($i=1,2$), it follows from   (\ref{5:vi})   that
\begin{equation}\label{5:bound}
  \sum_{i=1}^2\inte V_i(x)|u_{in}|^2dx= \sum_{i=1}^2\inte V_i(\epsilon_nx+\epsilon_ny_{\epsilon_n})|w_{in}|^2dx \leq \hat e(a^*-\beta,a^*-\beta,\beta)+\frac{1}{n}.
\end{equation}
We  deduce from  (\ref{5:R}) and Fatou's Lemma that
$\{\epsilon_ny_{\epsilon_n}\}$ is  bounded uniformly  in $\R^2$.

For any $\varphi(x)\in C_0^\infty(\R^2)$, define
\begin{equation*}\tilde \varphi(x)=\varphi\big(\frac{x-\epsilon_ny_{\epsilon_n}}{\epsilon_n}\big),\,\
j(\tau,\sigma)=\frac{1}{2}\inte\big|u_{1n}+\tau u_{1n}+\sigma\tilde \varphi\big|^2dx,
\end{equation*}
so that $j(\tau,\sigma)$ satisfies
$$j(0,0)=\frac{1}{2},\,\ \frac{\partial j(0,0)}{\partial \tau}=\inte|u_{1n}|^2dx=1\,\ \text{and}\,\ \frac{\partial j(0,0)}{\partial \sigma}=\inte u_{1n}\tilde \varphi dx.$$
Applying the implicit function theorem then gives that there exist a constant $\delta_n>0$ and a function $\tau(\sigma)\in C^1\big((-\delta_n,\delta_n), \R\big)$ such that
$$\tau(0)=0,\,\ \tau'(0)=-\inte u_{1n}\tilde \varphi dx, \,\ \text{and}\ \ j(\tau(\sigma),\sigma)=j(0,0)=\frac{1}{2}.$$
This implies that
$$\big(u_{1n}+\tau(\sigma) u_{1n}+\sigma\tilde \varphi, u_{2n}\big)\in \mathcal{M},\,\ \text{where}\,\ \sigma\in (-\delta_n,\delta_n).$$
We then obtain from (\ref{eqE2}) that
$$E_{a^*,a^*}(u_{1n}+\tau(\sigma) u_{1n}+\sigma\tilde \varphi, u_{2n})-E_{a^*,a^*}(u_{1n}, u_{2n})\geq-\frac{1}{n}\|(\tau(\sigma)u_{1n}+\sigma\tilde\varphi,0)\|_\mathcal{X}.$$
Setting $\sigma\to 0^+$ and $\sigma\to 0^-$, respectively, we thus have
\begin{equation}\label{5:de}\Big|\big\langle E_{a^*,a^*}'(u_{1n},u_{2n}), (\tau'(0)u_{1n}+\tilde \varphi, 0)\big\rangle\Big|\leq\frac{1}{n}\|\tau'(0)u_{1n}+\tilde \varphi\|_{\mathcal{H}_1}.
\end{equation}
By the definition of (\ref{5:scal}), direct calculations  yield that for $n\to\infty$,
\begin{equation}\label{5:de1}
\begin{split}
\frac{1}{2}\big\langle E_{a^*,a^*}'(u_{1n},u_{2n}), (\tilde \varphi, 0)\big\rangle&=\frac{1}{\epsilon_n}\inte \nabla w_{1n}\nabla \varphi dx+\epsilon_n \inte V_1(\epsilon_n x+\epsilon_ny_{\epsilon_n})w_{1n}\varphi dx\\
&-\frac{a^*}{\epsilon_n}\inte w_{1n}^3\varphi dx+\frac{\beta}{\epsilon_n}\inte (w_{1n}^2-w_{2n}^2)w_{1n}\varphi dx,
\end{split}
\end{equation}
and
\begin{equation}\label{5:de2}
\begin{split}
&\tau'(0)=-\inte u_{1n}\tilde\varphi dx=-\epsilon_n\inte w_{1n}\varphi dx\,,\quad \|\tau'(0)u_{1n}+\tilde \varphi\|_{\mathcal{H}_1}<C,\\
&\mu_{1n}:=\frac{1}{2}\langle E_{a^*,a^*}'(u_{1n},u_{2n}), (u_{1n}, 0)\rangle\sim -\frac{a^*}{2}\inte |u_{1n}|^4 dx=-\frac{a^*}{2}\epsilon_n^{-2}.
\end{split}\end{equation}
We then deduce from (\ref{5:de})-(\ref{5:de2}) that
\begin{equation}\label{5:de3}
\begin{split}
&\bigg|\inte \nabla w_{1n}\nabla \varphi dx+\epsilon_n^2 \inte V_1(\epsilon_n x+\epsilon_ny_{\epsilon_n})w_{1n}\varphi dx-\mu_{1n}\epsilon_n^2\inte w_{1n}\varphi dx
\\&-a^*\inte w_{1n}^3\varphi dx+\beta\inte (w_{1n}^2-w_{2n}^2)w_{1n}\varphi dx\bigg|\leq \frac{C\epsilon_n}{n}.
\end{split}\end{equation}
Using  (\ref{5:R}) and (\ref{5:L2}), we thus deduce from (\ref{5:de3}) that $w_{1n}\overset{n}{\rightharpoonup} w_1$  in $H^1(\R^2)$, where $w_1$ is a non-zero solution of
\begin{equation}\label{5:lim}
-\Delta w+\lambda^2w -a^*w^3=0 \,\ \text{in}\,\ \R^2,
\end{equation}
and $\lambda^2:=-\lim_{n\to\infty}\mu_{1n}\epsilon_n^2>0$. Similar to  the above argument, one can also show that $w_{2n}\overset{n}{\rightharpoonup} w_2$ in $H^1(\R^2)$, where $w_2$ is also a non-zero solution of (\ref{5:lim}). Further, we derive from (\ref{5:L2}) that $w_1(x)\equiv \pm w_2(x) $  a.e. in $ \R^2$.

%Since positive solutions of (\ref{5:lim}) must be radially symmetric (cf. \cite{Li}) and unique (cf. \cite{K}), we finally conclude that
%\begin{equation}\label{5:cons}
%w_1(x)\equiv w_2(x) =\frac{|\lam| }{a^*}Q(|\lam|x)\quad \mbox{in}\quad \R^2 ,\end{equation}
%and therefore, $\|w_{1}\|_2^2=\|w_{2}\|_2^2=1$.

We next show that
\begin{equation}\label{5:cons}
\|w_{i}\|_2^2=1,\,\ \mbox{where}\,\ i=1,\, 2.
\end{equation}
Indeed, since $\|w_{in}\|_2^2\equiv1$ and $w_i\not\equiv0$, we have $0<\|w_{i}\|_2^2\leq1$. On the other hand,   employing  (\ref{5:lim}) and the  Pohozaev identity  (\cite[Lemma 8.1.2]{C}), we derive that
$$\inte|w_i|^2dx=\frac{1}{\lambda^2}\inte |\nabla w_i|^2dx=\frac{a^*}{2\lambda^2}\inte|w_i|^4dx,\,\ \mbox{where}\,\ i=1,\, 2.$$
 We then use the Gagliardo-Nirenberg inequality (\ref{GNineq}) to deduce that
\begin{equation}
\frac{a^*}{2}\leq \frac{\inte|w_i|^2dx\inte |\nabla w_i|^2dx}{\inte|w_i|^4dx}=\frac{a^*}{2}\inte|w_i|^2dx,
\end{equation}
which yields that $\|w_{i}\|_2^2\geq1$ for $i=1$ and $2$, and therefore (\ref{5:cons}) follows.

We now obtain from  (\ref{5:cons})  that
\begin{equation}\label{5:strong}w_{in}\xrightarrow{n} w_i\,\ \text{strongly in}\,\ L^2(\R^2), \,\ i=1,\, 2.\end{equation}
Since $\{\epsilon_n y_{\epsilon_n}\}$ is  bounded uniformly in $\R^2$, there exists a subsequence, still denoted by $\{\epsilon_n\}$, of $\{\epsilon_n\}$ such that $\epsilon_n y_{\epsilon_n}\xrightarrow{n}z_0\in \R^2$.
Using Fatou's lemma, we thus obtain from (\ref{5:bound}) and (\ref{5:strong}) that
\begin{equation*}
\begin{split}
\hat{e}(a^*-\beta,a^*-\beta,\beta)&\geq \sum_{i=1}^2\inte\lim_{n\to\infty} V_i(\epsilon_nx+\epsilon_ny_{\epsilon_n})|w_{in}|^2dx\\
&= \sum_{i=1}^2\inte V_i(z_0)|w_{i}|^2dx=V_1(z_0)+V_2(z_0),
\end{split}
\end{equation*}
which however contradicts the assumption that $\hat{e}(a^*-\beta,a^*-\beta,\beta)<\inf_{x\in \R^2}\big(V_1(x)+V_2(x)\big)$, and the proof of Theorem \ref{thm03}(ii) is therefore done.\qed

%\begin{rem}
%By the norm preservation we further conclude
%that $w_{in}$ converges to $w_i$ strongly
%in $L^p(\mathbb{R}^2)$ for any $2\leq
%p<\infty$ because of $H^1(\mathbb{R}^2)$ boundedness. Also, since
%$w_{a_i}$ and $w_0$ satisfy (\ref{5:de3}) and (\ref{5:lim}),
%respectively, a simple analysis shows that $w_{a_i}$ converges to
%$w_0$ strongly in $H^1(\mathbb{R}^2)$. Moreover, we see that $w_i$ is a function such that  the equality holds for the Gagliardo-Nirenberg inequality (\ref{GNineq}), then by the uniqueness for this equality sated in \cite{W}, we see that
%$$w_i(x)=\pm \frac{\lambda}{\sqrt {a^*}}Q(\lambda|x-x_0|)\quad \text{for some }\quad x_0\in\R^2\,.  $$
%\end{rem}

\section{Uniqueness of nonnegative minimizers}

 The aim of this section is to  prove Theorem \ref{thm1.3} on the uniqueness of non-negative minimizers for (\ref{eq1.4}) by applying the implicit function theorem. Based on  the contracting map, a similar result for single GP energy functionals was also proved in \cite{Schnee,M}, however this methods  seems not applicable to our problem.

Let $\lambda_{i1} $  be the first eigenvalue of $-\Delta+V_i(x)$ in $\h_i$, $i.e,$
\begin{equation}\label{u:mu1}
\lambda_{i1}=\inf\Big\{\inte \big(|\nabla u|^2+V_i(x)\Big)u^2dx:\  u\in\h_i \,\text{ and }\,\inte |u|^2dx=1\Big\},\,\
 i=1,\, 2.\end{equation}
Applying
Lemma \ref{2:lem1} (see also \cite{RS}), one can deduce that $\lambda_{i1}$ is simple and can be attained by a positive normalized function $\phi_{i1}\in\h_i$, which is called the first eigenfunction of $-\Delta+V_i(x)$ in $\h_i$, $i=1,\, 2$.
Define
 $$Z_i= \text{span}\{\phi_{i1}\}^\perp=\Big\{ u\in\h_i:\inte u \phi_{i1}dx=0\Big\}, \,\ i=1,2,$$
so that
 \begin{equation}\label{u:sep}\h_i= \text{span}\{\phi_{i1}\}\oplus Z_i,\,\ i=1,2.\end{equation}
We now recall the following properties.

\begin{prop}\label{u:lem2}\textbf{\cite[Lemma 4.2]{GZZ}}  If $V_i(x)$ satisfies   (\ref{1.12}) for $i=1,\,2$, then
\begin{enumerate}
   \item [\rm(i)]  ker$\big(-\Delta+V_i(x)-\lambda_{i1}\big)=\text{span}\{\phi_{i1}\}$,
   \item [\rm(ii)] $\phi_{i1}\not\in \big(-\Delta+V_i(x)-\lambda_{i1}\big)Z_i$,
   \item [\rm(iii)] Im $\big(-\Delta+V_i(x)-\lambda_{i1}\big)$=$\big(-\Delta+V_i(x)-\lambda_{i1}\big)Z_i$ is closed in $\h^*_i$,
   \item[\rm(iv)] codim Im $\big(-\Delta+V_i(x)-\lambda_{i1}\big)$=1,
 \end{enumerate}
 where $\h^*_i$ denotes the dual space of $\h_i$ for $i=1,\, 2$.
\end{prop}
Define now the  $C^1$ functional  $F_i: \mathcal{X}\times \R^3\mapsto \h_i^*$ by
 \begin{equation}\label{u:eq1}F_i(u_1,u_2, \mu_i, b_i,\beta)=\big(-\Delta +V_i(x)-\mu_i\big)u_i-b_i u_i^3-\beta u_j^2u_i, \ \  i=1,2,\end{equation}
 where $j\not=i$ and $j=1,2$.
We then have the following lemma.

\begin{lem}\label{u:le3}
Let $F_i$ be defined by (\ref{u:eq1}), where $i=1,\, 2$. Then there exist $\delta>0$ and  a  unique function $(u_i(b_1,b_2,\beta), \mu_i(b_1,b_2,\beta))\in C^1\big(B_\delta(\vec0); B_\delta(\phi_{i1}, \lambda_{i1})\big)$, where $i=1,2$,
 such that
\begin{equation} \label{u:F}
\begin{cases}
\mu_i(\vec0)=\lambda_{i1},\ \,  u_i(\vec0)=\phi_{i1},\ i=1,\,2;\\
F_i\big(u_1(b_1,b_2,\beta),u_2(b_1,b_2,\beta), \mu_i(b_1,b_2,\beta), b_i,\beta\big)=0,\ \, i=1,\,2;\\
\|u_1(b_1,b_2,\beta)\|_2^2=\|u_2(b_1,b_2,\beta)\|_2^2=1.
\end{cases}
\end{equation}
\end{lem}

 \noindent{\bf Proof.} For $i=1,2$, define $g_i: (Z_1\times \R)\times (Z_2\times\R)\times\R^4\mapsto \h_i^*$ by
 $$g_i\big((z_1,\mu_1), (z_2,\mu_2),s_1,s_2,b_i,\beta\big):=F_i\big((1+s_1)\phi_{11}+z_1,(1+s_2)\phi_{21}+z_2,\mu_i, b_i,\beta\big).$$
Then  $g_i\in C^1((Z_1\times \R)\times (Z_2\times\R)\times\R^4; \h_i^*)$ and
\begin{equation}
 \begin{split}
 g_i((0,\lambda_{11}),(0,\lambda_{21}),\vec0)&=F_i(\phi_{11}, \phi_{21},\lambda_{i1},\vec0)=0,\\
D_{s_i}g_i((0,\lambda_{11}),(0,\lambda_{21}),\vec0)&=D_{u_i}F_i(\phi_{11}, \phi_{21},\lambda_{i1},\vec0)\phi_{i1} \\ &=\big(-\Delta+V_i(x)-\lambda_{i1}\big)\phi_{i1}=0,\\
D_{s_j}g_i((0,\lambda_{11}),(0,\lambda_{21}),\vec0)&=D_{u_j}F_i(\phi_{11}, \phi_{21},\lambda_{i1},\vec0)\phi_{j1} =0,
 \end{split}
 \end{equation}
where $j\not=i,i,j=1,\,2$. Moreover, for any $(\hat z_i, \hat \mu_i)\in Z_i\times \R$, we have
 \begin{equation}\label{eq4.13}
 \begin{split}
 &\big\langle D_{(z_i,\mu_i)}g_i((0,\lambda_{11}),(0,\lambda_{21}),\vec0),(\hat z_i, \hat \mu_i)\big\rangle\\
 =&D_{u_i}F_i(\phi_{11}, \phi_{21},\lambda_{i1},\vec0)\hat z_i +D_{\mu_i}F_i(\phi_{11}, \phi_{21},\lambda_{i1},\vec0)\hat \mu_i\\
 =&\big(-\Delta+V(x)-\lambda_{i1}\big)\hat z_i-\hat\mu_i \phi_{i1}\in \h_i^*,
 \end{split}
 \end{equation}
and
\begin{equation}
 \begin{split}
 &\big\langle D_{(z_j,\mu_j)}g_i((0,\lambda_{11}),(0,\lambda_{21}),\vec0),(\hat z_j, \hat \mu_j)\big\rangle\\
 =&D_{u_j}F_i(\phi_{11}, \phi_{21},\lambda_{i1},\vec0)\hat z_j +D_{\mu_j}F_i(\phi_{11}, \phi_{21},\lambda_{i1},\vec0)\hat \mu_j=0.
 \end{split}
 \end{equation}
It then follows from  (\ref{eq4.13}) and Proposition \ref{u:lem2} that for $i=1,2,$
\begin{equation}\label{4.15}D_{(z_i,\mu_i)}g_i((0,\lambda_{11}),(0,\lambda_{21}),\vec0): Z_i\times \R\mapsto \h_i^* \text{ is an isomorphism.}
\end{equation}

We next define $G: (Z_1\times \R)\times (Z_2\times\R)\times\R^5\mapsto \h_1^*\times \h_2^*\times \R^2$ by
\begin{equation*}
G\big((z_1,\mu_1), (z_2,\mu_2),s_1,s_2,(b_1,b_2,\beta)\big):=\left(\begin{array}{c}
g_1\big((z_1,\mu_1), (z_2,\mu_2),s_1,s_2,b_1,\beta\big) \\
 g_2\big((z_1,\mu_1), (z_2,\mu_2),s_1,s_2,b_2,\beta\big) \\
 \|(1+s_1)\phi_{11}+z_1\|_2^2-1\\
 \|(1+s_2)\phi_{21}+z_2\|_2^2-1
\end{array}
 \right).
\end{equation*}
Setting $h_i(z_i,s_i)=\|(1+s_i)\phi_{i1}+z_i\|_2^2-1$, we then have
\begin{equation}
\begin{split}
&D_{((z_1,\mu_1),(z_2,\mu_2),s_1,s_2)}G\big((0,\lambda_{11}), (0,\lambda_{21}),0,0,\vec 0\,\big)\\
=&
\left(
  \begin{array}{cccc}
    D_{(z_1,\mu_1)}g_1 & D_{(z_2,\mu_2)}g_1 & D_{s_1}g_1 & D_{s_2}g_1 \\
    D_{(z_1,\mu_1)}g_2 & D_{(z_2,\mu_2)}g_2 & D_{s_1}g_2 & D_{s_2}g_2 \\
    D_{(z_1,\mu_1)}h_1& D_{(z_2,\mu_2)}h_1 & D_{s_1}h_1&  D_{s_2}h_1\\
   D_{(z_1,\mu_1)}h_2& D_{(z_2,\mu_2)}h_2& D_{s_1}h_2&  D_{s_2}h_2 \\
  \end{array}
\right)\\
=&\left(
  \begin{array}{cccc}
    (-\Delta+V_1(x)-\lambda_{11}, -\phi_{11}) &0 & 0 & 0 \\
    0 & (-\Delta+V_2(x)-\lambda_{21}, -\phi_{21}) & 0 & 0 \\
    0 & 0 & 2& 0 \\
    0 & 0 & 0& 2 \\
  \end{array}
\right).
\end{split}
\end{equation}
We then derive from (\ref{4.15})  that
\begin{equation*}
D_{((z_1,\mu_1),(z_2,\mu_2),s_1,s_2)}G\big((0,\lambda_{11}), (0,\lambda_{21}),0,0,\vec 0\big):
(Z_1\times \R)\times (Z_2\times\R)\times\R^5
\mapsto \h_1^*\times \h_2^*\times\R^2
\end{equation*}is an isomorphism.
Therefore, it follows from  the implicit function theorem that there exist $\delta>0$ and a unique function $(z_i(b_1,b_2,\beta),\mu_i(b_1,b_2,\beta),s_i(b_1,b_2,\beta))\in C^1(B_{\delta}(\vec 0);
B_{\delta}(0,\lambda_{i1},0))$, where $i=1,\,2$, such that
  \begin{equation}\label{u:G}
 \begin{cases}
 G\big((z_1,\mu_1), (z_2,\mu_2),s_1,s_2,(b_1,b_2,\beta)\big)=G\big((0,\lambda_{11}), (0,\lambda_{21}),0,0,\vec0\big)=\vec 0,\\
 z_i(\vec0)=0,\,\ \mu_i(\vec0)=\lambda_{i1},\,\ s_i(\vec0)= 0,\,\ i=1,\,2.\\
 \end{cases}
 \end{equation}
By setting
$$u_i(b_1,b_2,\beta)=(1+s_i(b_1,b_2,\beta))\phi_{i1}+z_i(b_1,b_2,\beta)\,,\ (b_1,b_2,\beta)\in B_{\delta}(\vec0),\ i=1,\,2,$$
we then obtain from (\ref{u:G})  that there exists  a unique function $(u_i(b_1,b_2,\beta),\mu_i(b_1,b_2,\beta))\in C^1(B_{\delta}(\vec 0);B_{\delta}(\phi_{i1},\lambda_{i1}))$, where $i=1,2$,  such that
\begin{equation}
u_i(\vec0)=(1+s_i(\vec0))\phi_{i1}+z_i(\vec0)=\phi_{i1},\,\ \mu_i(\vec0)=\lambda_{i1},
\end{equation}
and
\begin{equation}
\left(\begin{array}{c}
 F_1\big(u_1(b_1,b_2,\beta),u_2(b_1,b_2,\beta),\mu_1(b_1,b_2,\beta), b_1,\beta\big) \\
  F_2\big(u_1(b_1,b_2,\beta),u_2(b_1,b_2,\beta),\mu_2(b_1,b_2,\beta), b_2,\beta\big) \\
   \|u_1(b_1,b_2,\beta)\|_2^2-1\\
     \|u_1(b_1,b_2,\beta)\|_2^2-1
      \end{array}
      \right)=\vec 0,
\end{equation}
and therefore (\ref{u:F}) holds. This completes the proof of Lemma \ref{u:le3}.
\qed

%\subsection{Existence and uniqueness of minimizers}

%%In fact, the existence of the minimizers of problem was proved in \cite{GS}, and we only  need to show the uniqueness of non-negative  minimizers.

In the following we use Lemma \ref{u:le3} to derive the uniqueness of nonnegative minimizers for sufficiently small $|(b_1,b_2,\beta)|$.

\noindent\textbf{Proof of Theorem \ref{thm1.3}.}
It follows from Theorem \ref{thm1} that $\hat e(b_1,b_2,\beta)$ admits at least one minimizer if $0<b_1<a^*$, $0<b_2<a^*$ and  $\beta<\sqrt {(a^*-b_1)(a^*-b_2)}$. We first claim that $\hat e(\cdot)$ is a continuous function of $(b_1,b_2,\beta)$ on the interval $I:=(-\frac{a^*}{2},\frac{a^*}{2})\times(-\frac{a^*}{2},\frac{a^*}{2})\times(-\frac{a^*}{4},\frac{a^*}{4})$. Indeed, for any $(b_1,b_2,\beta)\in I$, let  $(u_1,u_2)$ be any nonnegative minimizer of $\hat e(b_1,b_2,\beta)$. It then follows from  (\ref{GNineq}) and   Cauchy's inequality that
\begin{equation*}
\begin{split}
\hat e(b_1,b_2,\beta)=E_{b_1,b_2,\beta}(u_1,u_2)
\geq \sum_{i=1}^2\big(\frac{a^*-|b_i|}{2}-\frac{|\beta|}{2}\big)\inte |u_i|^4dx\geq \frac{a^*}{8}\sum_{i=1}^2\inte |u_i|^4dx,
\end{split}
\end{equation*}
which implies that
\begin{equation}\label{4.19}\text{the $L^4$-norm of  minimizers for $\hat e(\cdot)$ is bounded uniformly on the interval $I$}.\end{equation}
For any $(b_{i1},b_{i2},\beta_i)\in I$, we now denote $(u_{i1},u_{i2})$ to be the corresponding nonnegative minimizer of $\hat e(b_{i1},b_{i2},\beta_i)$, where $i=1,\,2$. Then
\begin{equation}\label{4.20}\begin{split}
\hat e(b_{11},b_{12},\beta_1)&=E_{b_{11},b_{12},\beta_1}(u_{11},u_{12}) \\
&=E_{b_{21},b_{22},\beta_2}(u_{11},u_{12})
+\frac{b_{21}-b_{11}}{2}\inte |u_{11}|^4dx \\  \quad &\quad +\frac{b_{22}-b_{12}}{2}\inte|u_{12}|^4dx+(\beta_2-\beta_1)\inte|u_{11}|^2|u_{12}|^2dx \\
\quad &\geq \hat e(b_{21},b_{22},\beta_2)+O(|(b_{21},b_{22},\beta_2)-(b_{11},b_{12},\beta_1)|).
\end{split}\end{equation}
Similarly, we also have
\begin{align}
&\hat e(b_{21},b_{22},\beta_2)\geq \hat e(b_{11},b_{12},\beta_1)+O(|(b_{21},b_{22},\beta_2)-(b_{11},b_{12},\beta_1)|).\label{4.21}
\end{align}
We then derive from (\ref{4.20}) and (\ref{4.21}) that
$$\lim_{(b_{11},b_{12},\beta_1)\to(b_{21},b_{22},\beta_2)}\hat e(b_{11},b_{12},\beta_1)=\hat e(b_{21},b_{22},\beta_2),$$
which implies that $\hat e(\cdot)$ is continuous  on the interval $I$, and the above claim is therefore proved.

Let $(u_{b_1,\beta},u_{b_2,\beta})$ be a non-negative minimizer of $\hat e(b_1,b_2,\beta)$ with $(b_1,b_2,\beta)\in I$.
We then deduce that $(u_{b_1,\beta},u_{b_2,\beta})$ satisfies the  Euler-Lagrange system
\begin{equation}\label{4.24}
\begin{cases}
 (-\Delta+V_1(x)-\mu_{b_1,\beta}) u_{b_1,\beta}-b_1
u_{b_1,\beta}^3-\beta u_{b_2,\beta}^2u_{b_1,\beta}=0\,\ \mbox{in}\,\ \R^2,\\[2mm]
(-\Delta+V_2(x)-\mu_{b_2,\beta}) u_{b_2,\beta}-b_2
u_{b_2,\beta}^3-\beta u_{b_1,\beta}^2u_{b_2,\beta}=0\,\ \mbox{in}\,\ \R^2,
\end{cases}
\end{equation}
i.e, \begin{equation}\label{4.24s}
F_i(u_{b_1,\beta},u_{b_2,\beta}, \mu_{b_i,\beta}, b_i,\beta)=0, \,\ i=1,2,
\end{equation}
where $(\mu_{b_1,\beta}, \mu_{b_2,\beta})\in\R ^2$ is a Lagrange multiplier.
We then derive from (\ref{4.19}) and the above claim that
\begin{equation}\begin{split}
&E_{0,0,0}(u_{b_1,\beta},u_{b_2,\beta})\\
=&E_{b_1,b_2,\beta}(u_{b_1,\beta},u_{b_2,\beta})+\sum_{i=1}^2\frac{b_i}{2}\inte|u_{b_i,\beta}|^4dx+
\beta\inte|u_{b_1,\beta}|^2|u_{b_2,\beta}|^2dx \\
=&\hat e(b_1,b_2,\beta)+O(|(b_1,b_2,\beta)|)\to \hat e(0,0,0)\,\ \text{ as }\,\ (b_1,b_2,\beta)\to 0.
\end{split}\label{4.23}\end{equation}

On the other hand, one can check easily that
$\hat e(0,0,0)=\lambda_{11}+\lambda_{21},$
 and $(\phi_{11}, \phi_{21})$ is the unique non-negative minimizer of $\hat e(0,0,0)$,
where $(\lambda_{i1},\phi_{i1})$ is the first eigenpair of $-\Delta+V_i(x)$ in $\h_i$, $i=1,\,2$. We then deduce from Lemma \ref{2:lem1} that
\begin{equation}\label{4.26}
u_{b_i,\beta}\to \phi_{i1}\ \text{ in } \ \, \h_i\,\ \text{ as } \ \, (b_1,b_2,\beta)\to(0,0,0),\,\ i=1,\, 2.
\end{equation}
Moreover, by (\ref{4.24}) and (\ref{4.26})  we have
\begin{align}
\mu_{b_i,\beta}&=\inte |\nabla u_{b_i,\beta}|^2+V_i(x)|u_{b_i,\beta}|^2dx-b_i\inte|u_{b_i,\beta}|^4dx-\beta\inte|u_{b_1,\beta}|^2|u_{b_2,\beta}|^4dx\nonumber\\
&\to \lambda_{i1}\quad \text{ as } \ (b_1,b_2,\beta)\to(0,0,0),\quad i=1,\,2.\label{4.27}
\end{align}
Applying (\ref{4.26}) and (\ref{4.27}), we then derive that there exists a small constant $\delta_1>0$ such that
\begin{equation}\label{4.28}\|u_{b_i,\beta}-\phi_{i1}\|_{\h_i}<\delta_1\,\ \text{ and } \,\ |\mu_{b_i,\beta}-\lambda_{i1}|<\delta_1 \,\ \text{ if }\, \ (b_1,b_2,\beta)\in B_{\delta_1}(\vec 0),\ \  i=1,\,2.\end{equation}
We thus conclude from (\ref{4.24s}) and Lemma \ref{u:le3} that
$$\mu_{b_i,\beta}=\mu_i(b_1,b_2,\beta), \ u_{b_i,\beta}=u_i(b_1,b_2,\beta), \,\ \text{ if }\, \ (b_1,b_2,\beta)\in B_{\delta_1}(\vec0),\, \  i=1,\,2.$$
This therefore implies that for sufficiently small $|(b_1,b_2,\beta)|$, $(u_1(b_1,b_2,\beta),\, u_2(b_1,b_2,\beta))$ is a unique non-negative minimizer of $\hat e(b_1,b_2,\beta)$,  and we are done.
\qed

\section{Optimal Estimates }\label{Sec5}

 To simplify the notations and the proof,  we denote $a_i=b_i+\beta>0$ for any fixed $0<\beta <a^*$, where $i=1,\, 2$. It is then equivalent to rewriting the functional (\ref{f1}) as
\begin{equation}\label{f}
\begin{split}
  E_{a_1,a_2}(u_1,u_2):=&\sum_{i=1}^2\int_{\R ^2} \Big(|\nabla
  u_i|^2+V_i(x)|u_i|^2-\frac{a_i}{2}|u_i|^4\Big)dx\\
  &+\frac{\beta}{2}\inte \big(|u_1|^2-|u_2|^2\big)^2dx \,,\quad\mbox{where}\quad  (u_1,u_2)\in\mathcal{X}.
\end{split}\end{equation}
Also, the minimization problem (\ref{eq1.4}) is then equivalent  to the following one:
\begin{equation}\label{eq1.3}
e(a_1,a_2):=\inf_{\{(u_1,u_2)\in \mathcal{M}\}} E_{a_1,a_2}(u_1,u_2),
\end{equation}
where  $\mathcal{M}$ is defined by (\ref{norm}).

To prove Theorem \ref{thm1.4}, we first need to establish  the following crucial optimal energy estimates  as $(a_1,a_2)\nearrow(a^*, a^*)$.

\begin{prop}\label{3:thm2} Suppose $0<\beta <a^*$ and  $V_i(x) $ satisfies (\ref{1:V}) and (\ref{1:V1}), where $i=1,\,2$. Then there exist two positive
constants $C_1$ and $ C_2$, independent of $a_1$ and $a_2$, such that
\begin{equation}\label{energy}
C_1\Big(a^*-\frac{a_{1}+a_{2}}{2}\Big)^\frac{p_0}{p_0+2}\leq e(a_1,a_2)\leq C_2
\Big(a^*-\frac{a_{1}+a_{2}}{2}\Big)^\frac{p_0}{p_0+2}\quad
\mbox{as}\quad
 (a_1,a_2)\nearrow (a^*, a^*),
\end{equation}
where  $p_0>0$ is defined by (\ref{1:def}), and $e(a_1,a_2)$ is defined by (\ref{eq1.3}).
Moreover, if $(u_{a_1},u_{a_2})$ is a non-negative minimizer of
$e(a_1,a_2)$, then there exist two positive
constants $C_3$ and $ C_4$, independent of $a_1$ and $a_2$, such that
\begin{equation}\label{1:u4}
C_3\Big(a^*-\frac{a_{1}+a_{2}}{2}\Big)^{-\frac{2}{p_0+2}}\leq \inte|u_{a_i}|^4dx\leq C_4
\Big(a^*-\frac{a_{1}+a_{2}}{2}\Big)^{-\frac{2}{p_0+2}}
\end{equation}
as $(a_1,a_2)\nearrow (a^*,a^*)$.
\end{prop}

We remark that even though the  upper bound of (\ref{energy}) can be proved similarly to that  of Lemma 3 in \cite{GS}, the arguments of \cite{GS} do  not give the  lower bound of (\ref{energy}). For this reason, as discussed below, we need employ a little more involved analysis  to address the optimal lower bound of (\ref{energy}).

%\subsection{Optimal energy estimates}
In what follows, we focus on the proof of Proposition \ref{3:thm2}.
For any fixed $0<\beta <a^*$, denote $(u_{a_1},u_{a_2})$ to be a non-negative minimizer of (\ref{eq1.3}). We start with the following energy estimates of $e(a_1,a_2)$.

\begin{lem}\label{lem3.1}
Under the assumptions of Proposition \ref{3:thm2},  there exists a constant $C>0$, independent of $a_1$ and $a_2$, such that
\begin{equation}\label{3:est}
e_1(a_1)+e_2(a_2)+\frac{\beta}{2}\int_{\R^2}\big(|u_{a_1}|^2-|u_{a_2}|^2\big)^2dx\leq
e(a_1,a_2)\leq C
\Big(a^*-\frac{a_1+a_2}{2}\Big)^\frac{p_0}{p_0+2}
\end{equation}
as $(a_1, a_2)\nearrow (a^*,a^*)$, where $e_i(\cdot)$ is given by (\ref{1:two}) for $i=1,\,2$.
\end{lem}

\noindent{\bf Proof.} Since $(u_{a_1},u_{a_2})$ is a non-negative minimizer of
(\ref{eq1.3}), we note  from (\ref{single}) that
\begin{equation}\label{3:equ}\nonumber
E_{a_1,a_2}(u,v)=E_{a_1}^1(u)+E_{a_2}^2(v)+\frac{\beta}{2}
\int_{\R^2}\big(|u|^2-|v|^2\big)^2dx \quad\text{for all}~ (u,v)\in
\mathcal{X},\end{equation}
where $E_{a_i}^i(\cdot)$ is defined by (\ref{single}) for $i=1$ and $2$.
This relation then implies that
\begin{equation}\label{3:left}\nonumber
e(a_1,a_2)\geq e_1(a_1)+e_2(a_2)+\frac{\beta}{2}
\int_{\R^2}\big(|u_{a_1}|^2-|u_{a_2}|^2)^2dx,
\end{equation}
which gives the lower bound of (\ref{3:est}).

Stimulated by Lemma 3 in \cite{GS}, we next prove the upper bound of (\ref{3:est}) as follows. Without loss of
generality, we may assume $p_0=\bar p_1=\min\{p_{11},p_{21}\}>0$ and $p_{11}\leq p_{21}$, where  $p_0$ and $\bar p_i$ are defined
 by (\ref{1:def}).
We proceed similarly to the proof of Lemma A2 in the Appendix, and use the
trial function $(\phi,\phi)$ for $\phi$ satisfying (\ref{2:trial}) with $\bar x_0 =x_{11}$.
Choose  $R>0$ small enough  that
$$V_i(x)\leq
C|x-x_{11}|^{p_{i1}}\quad \text{for}\quad |x-x_{11}|\leq R,\quad
i=1,\,2.$$
By the exponential decay of $Q(x)$, we have
\begin{equation*}
\begin{split}
\int_{\R^2}V_i(x)\phi^2(x)dx\leq
C\tau^{-p_{i1}}\inte|x|^{p_{i1}}Q^2(x)dx\leq C\tau^{-p_{i1}}\quad
 {\rm as}\quad \tau \to \infty\,, \quad i=1,\, 2.
\end{split}
\end{equation*}
This inequality and (\ref{2:limt2}) then imply that
\begin{equation*}
\begin{split}
E_{a_1,a_2}(\phi,\phi)&\leq\frac{2}{a^*}\Big(a^*-\frac{a_1+a_2}{2}\Big)
\tau^2+C\big(\tau^{-p_{11}}+\tau^{-p_{21}}\big).
\end{split}
\end{equation*}
Setting
$\tau=\Big(a^*-\frac{a_1+a_2}{2}\Big)^\frac{-1}{p_0+2}$ and using
  $p_0=p_{11}\leq p_{21}$, we derive that
$$e(a_1,a_2)\leq E_{a_1,a_2}(\phi,\phi)\leq C
\Big(a^*-\frac{a_1+a_2}{2}\Big)^\frac{p_0}{p_0+2}\,,$$ which therefore gives
the upper bound of (\ref{3:est}).
\qed \\

By Lemma \ref{lem3.1}, for $i=1,\,2$, we have
\begin{equation} e_i(a_i)\leq E_{a_i}^i(u_{a_i})\leq C
\Big(a^*-\frac{a_1+a_2}{2}\Big)^\frac{p_0}{p_0+2}
 \label{3:a_1}
\end{equation}
as $(a_1, a_2)\nearrow (a^*,a^*)$, where $e_i(\cdot)$ is defined by  (\ref{1:two}).
On the other hand, it is proved in \cite[Lemma 3]{GS} that for
\begin{equation}\label{1:pi}p_i:=\max\{p_{ij},j=1,\ldots,
n_i\}>0\,,\quad i=1,\,2\,,\end{equation}  there exists two positive constants $m$ and $M$, independent of $a_1$ and $a_2$, such that
\begin{equation}\label{1:ei}
m (a^*-a_i)^\frac{p_i}{p_i+2}\leq e_i(a_i)\leq M
(a^*-a_i)^\frac{p_i}{p_i+2}\ \, \text{for}\ \, 0\leq a_i\leq
a^*\,,\ \mbox{and}\ \, i=1,\,2.
\end{equation}
Applying (\ref{1:ei}) and Lemma \ref{lem3.1}, we next derive the following  $L^4(\R^2)-$estimates of minimizers.

\begin{lem}\label{le3.2}
Under the assumptions of   Proposition \ref{3:thm2}, we have
{\small\begin{equation}\label{3:L2}
C\Big(a^*-\frac{a_1+a_2}{2}\Big)^{-\frac{ 2}{p_0+2}\frac{p_0}{p_i}}\leq
\int_{\R^2}|u_{a_i}|^4dx\leq
\frac{1}{C}\Big(a^*-\frac{a_1+a_2}{2}\Big)^{-\frac{ 2}{p_0+2}} \text{ as }(a_1, a_2)\nearrow (a^*,a^*),
\end{equation}}
and
\begin{equation}\label{3:L0}
\lim_{(a_1,a_2)\nearrow
(a^*,a^*)}\frac{\int_{\R^2}|u_{a_1}|^4dx}{\int_{\R^2}|u_{a_2}|^4dx}=1,
\end{equation}
 where $p_i\ge p_0$, and $p_i>0$ is given by (\ref{1:pi}) for $i=1 $ and $ 2$.
\end{lem}

\noindent{\bf Proof.}
We first prove the lower bound of  (\ref{3:L2}). Pick any $0<b<a_i<a^*$ ($i=1,\,2$), and observe that
\begin{equation*}
e_i(b)\leq
E^i_{a_i}(u_{a_i})+\frac{a_i-b}{2}\int_{\mathbb{R}^2}|u_{a_i}|^4dx,\quad
i=1,\, 2.
\end{equation*}
It then follows from (\ref{3:a_1}) and (\ref{1:ei}) that
\begin{equation}\label{eq2.7}
\begin{split}
\frac{1}{2}\int_{\mathbb{R}^2}|u_{a_i}|^4dx&\geq
\frac{e_i(b)-C\big(a^*-\frac{a_1+a_2}{2}\big)^\frac{p_0}{p_0+2}}{a_i-b}\\
&\geq
\frac{m(a^*-b)^\frac{p_i}{p_i+2}-C\big(a^*-\frac{a_1+a_2}{2}\big)^\frac{p_0}{p_0+2}}{a_i-b}.
\end{split}\end{equation}
Take
$b=a_i-C_1\big(a^*-\frac{a_1+a_2}{2}\big)^\frac{p_0(p_i+2)}{p_i(p_0+2)}$,
where $C_1>0$ is so large that
$mC_1^\frac{p_0(p_i+2)}{p_i(p_0+2)}>2C$. We then derive from
(\ref{eq2.7}) that
\begin{equation}\label{3:low}
\int_{\mathbb{R}^2}|u_{a_i}|^4dx\geq
C_2\big(a^*-\frac{a_1+a_2}{2}\big)^{-\frac{2}{p_0+2}\frac{p_0}{p_i}},\quad
i=1,\,2,
\end{equation}
which therefore implies the lower bound of (\ref{3:L2}).

%Without loss of
%generality, we may assume that
%\begin{equation}\label{3:ui}\int_{\R^2}|u_{a_1}|^4dx\leq\int_{\R^2}|u_{a_2}|^4dx.
%\end{equation}

We next  prove the upper bound of (\ref{3:L2}). One can note from  (\ref{3:a_1}) that
for $i=1,\,2$,
\begin{equation}\label{3:L1}
 E_{a_i}(u_{a_i})\leq C\Big(a^*-\frac{a_1+a_2}{2}\Big)^\frac{p_0}{p_0+2}\quad \mbox{as}
 \ \, (a_1,a_2)\nearrow (a^*,a^*).
\end{equation}
Without loss of generality, we may assume that ~$a_1\leq a_2\leq a^*$ and $(a_1,a_2)\not=(a^*,a^*)$. By  (\ref{GNineq}), we then have
\begin{equation*}
E_{a_1}^1(u_{a_1})\geq\frac{a^*-a_1}{2}\int_{\mathbb{R}^2}|u_{a_1}|^4dx\geq\frac{1}{2}
\Big(a^*-\frac{a_1+a_2}{2}\Big)\int_{\mathbb{R}^2}|u_{a_1}|^4dx.
\end{equation*}
It thus follows from (\ref{3:L1}) that  the upper bound of (\ref{3:L2}) holds for $u_{a_1}$. Similarly, the upper bound of (\ref{3:L2})  holds also for $u_{a_2}$ if (\ref{3:L0}) is true, and then the proof is done.

Now we come to prove (\ref{3:L0}). Recall from Lemma \ref{lem3.1} that
$$\int_{\R^2}\big(|u_{a_1}|^2-|u_{a_2}|^2\big)^2dx  \leq C
\Big(a^*-\frac{a_1+a_2}{2}\Big)^\frac{p_0}{p_0+2} \quad \text{as}\quad (a_1,a_2)\nearrow (a^*,a^*)\,,$$
which implies that
\begin{equation}\label{3:L3}\begin{split}
&\bigg|\int_{\R^2}|u_{a_1}|^4dx-
\int_{\R^2}|u_{a_2}|^4dx\bigg|\\
&\leq\Big(\int_{\R^2}\big(|u_{a_1}|^2-|u_{a_2}|^2\big)^2dx\Big)^
\frac{1}{2}\Big(\int_{\R^2}\big(|u_{a_1}|^2+|u_{a_2}|^2\big)^2dx\Big)^\frac{1}{2}\\
&\leq C
\Big(a^*-\frac{a_1+a_2}{2}\Big)^\frac{p_0}{2(p_0+2)}
\bigg[\Big(\int_{\R^2}|u_{a_1}|^4dx\Big)^\frac{1}{2}+\Big(\int_{\R^2}|u_{a_2}|
^4dx\Big)^\frac{1}{2}\bigg]
\end{split}\end{equation}
as $(a_1,a_2)\nearrow (a^*,a^*)$. Since it follows from (\ref{3:low}) that  $\int_{\R^2}|u_{a_i}|^4dx\to\infty$
as $(a_1,a_2)\nearrow (a^*,a^*)$, where $i=1,\,2$, we conclude (\ref{3:L0}) from the above estimate. \qed

%\subsection{Proof of Theorem \ref{3:thm2}}

%The approach of proving Theorem
%\ref{3:thm2} can be used effectively to estimate precisely GP energy
%functionals of other types, such as (\ref{f}) under general trapping
%potentials.

We next claim  that the upper estimates of (\ref{3:est})
and (\ref{3:L2}) are optimal. By Lemma \ref{lem3.1}, we see   that
\begin{equation}\label{eq2.8}
\sum_{i=1}^2\int_{\mathbb{R}^2}V_i(x)|u_{a_i}(x)|^2dx\leq
e(a_1,a_2)\leq C\big(a^*-\frac{a_1+a_2}{2}\big)^\frac{p_0}{p_0+2} \quad \text{as} \quad (a_1,a_2)\nearrow (a^*,a^*).
\end{equation}
Set
\begin{equation}\label{eq2.9a}\epsilon ^{-2}(a_1,a_2):=\int_{\mathbb{R}^2}|u_{a_1}(x)|^4dx,\quad \text{where}\ \ \epsilon(a_1,a_2)>0\,.
\end{equation}
It then yields from (\ref{3:L2}) that $\epsilon(a_1,a_2) \searrow 0$ as
$(a_1,a_2)\nearrow (a^*,a^*)$. Moreover, we deduce from (\ref{3:est}) and (\ref{3:L0})
that for $i=1,\,2$,
\begin{equation}\label{eq2.10}
\int_{\mathbb{R}^2}|
u_{a_i}|^4dx,\int_{\mathbb{R}^2}|\nabla u_{a_i}|^2dx\sim
\epsilon ^{-2}(a_1,a_2)\quad \text{as}\quad (a_1,a_2)\nearrow (a^*,a^*),
\end{equation}
where $f\sim g$ means that $f/g$ is bounded from below and above. In
view of above facts,  we next define the
$L^2(\mathbb{R}^2)$-normalized function
\begin{equation}\label{eq2.10a}
\tilde{w}_{a_i}(x):=\epsilon(a_1,a_2) u_{a_i}\big(\epsilon(a_1,a_2)
x\big),\ \, i=1,\,   2.
\end{equation}
It then follows from (\ref{3:L0}), (\ref{eq2.9a})  and
(\ref{eq2.10}) that for $i=1,\,2$,
 \begin{equation}\label{eq2.12}
\int_{\mathbb{R}^2}| \tilde{w}_{a_i}|^4dx=1\quad\text{and} \quad
m\leq \int_{\mathbb{R}^2}| \nabla \tilde{w}_{a_i}|^2dx\leq
\frac{1}{m}\quad \text{as}\quad (a_1,a_2)\nearrow (a^*,a^*),
\end{equation}
where $m>0$ is independent of $a_1$ and $a_2$. In the rest part of this
section, for convenience we use $\epsilon >0$ to denote $\epsilon(a_1,a_2)$ so that $\epsilon \searrow 0$ as $(a_1,a_2)\nearrow (a^*,a^*)$.

\begin{lem}\label{le3.3}
Under the assumptions of Proposition \ref{3:thm2}, we have
\begin{itemize}
  \item [\rm(i).]There exist a sequence $\{y_\epsilon\}$
as well as positive constants $R_0$ and $\eta_i$ such that for $i=1,\, 2$, the function
\begin{equation}\label{eq2.15}
 w_{a_i}(x)=\tilde{w}_{a_i}(x+y_\epsilon)=\epsilon
u_{a_i}(\epsilon x+\epsilon y_\epsilon)
\end{equation}
satisfies
\begin{equation}\label{eq2.16}
\lim_{\epsilon\to
0}\inf\int_{B_{R_0}(0)}|w_{a_i}|^2dx\geq\eta_i>0.
\end{equation}

\item [\rm(ii).] The following estimate holds
\begin{equation}\label{eq3.20}
  \lim_{(a_1,a_2)\nearrow(a^*,a^*)}{\rm dist}(\epsilon y_\epsilon, \Lambda)=0,
  \end{equation}
  where $\Lambda=\big\{x\in \R^2:\,
V_1(x)=V_2(x)=0\big\}$ is given by (\ref{1:com}).  Moreover, for
any sequence $\{(a_{1k},a_{2k})\}$ satisfying $(a_{1k},a_{2k})\nearrow (a^*,a^*)$ as $k\to\infty$,  there exists a   convergent
subsequence of $\{(a_{1k},a_{2k})\}$, still denoted by $\{(a_{1k},a_{2k})\}$, such that
\begin{equation}\label{eq3.21}
\epsilon_k y_{\epsilon_k}\xrightarrow{k} x_0\in \Lambda \text{ and }w_{a_{ik}}\overset{k}{\to} w_0\quad\text{strongly in } H^1(\R^2),\quad i=1,\ 2,\end{equation}
where $w_0$ satisfies
\begin{equation}\label{eq2.21}
w_0(x)=\frac{\lambda}{\|Q\|_2}Q(\lambda|x-y_0|) \quad \text{for
some} \  y_0\in\mathbb{R}^2\  \text{and}\ \, \lambda>0.\end{equation}
\end{itemize}
\end{lem}

\noindent \textbf{Proof.} \textbf{(i).} In  order to establish
(\ref{eq2.16}), in view of (\ref{eq2.15}) it suffices to prove that there exist $R_0>0$ and $\eta_i>0$ such that \begin{equation}\label{eq2.14}
\liminf_{\epsilon\to
0}\int_{B_{R_0}(y_\epsilon)}|\tilde{w}_{a_i}|^2dx\geq\eta_i>0\,,\quad \mbox{where}\quad
i=1,\,2.
\end{equation}
We first show that (\ref{eq2.14}) holds for $\tilde w_{a_1}$. Indeed, if it is false, then for any $R>0$, there exists
a subsequence $\{\tilde{w}_{a_{1k}}\}$, where $(a_{1k},a_{2k})\nearrow (a^*,a^*)$ as $k\to\infty$, such that
\begin{equation*}
\lim_{k\rightarrow\infty}\sup_{y\in
\mathbb{R}^2}\int_{B_{R}(y)}|\tilde{w}_{a_{1k}}|^2dx=0.
\end{equation*}
By Lemma I.1 in \cite{l2}, we then
deduce that $\tilde{w}_{a_{1k}}\xrightarrow{k}0$ in
$L^p(\mathbb{R}^2)$ for any $2<p<\infty$, which contradicts
(\ref{eq2.12}). Thus $\tilde w_{a_1}$ satisfies (\ref{eq2.14}) for a
sequence $\{y_\epsilon\}$, $R_0>0$ and $\eta _1>0$.

We next show that for the sequence $\{y_\epsilon\}$, $R_0>0$ and $\eta _1>0$
obtained above, (\ref{eq2.14}) holds also
for $\tilde w_{a_2}$ with some constant $\eta_2>0$. On the contrary, if
(\ref{eq2.14}) is false for $\tilde w_{a_2}$, then  there exists
a subsequence $\{\tilde{w}_{a_{2k}}\}$, where $(a_{1k},a_{2k})\nearrow (a^*,a^*)$ as $k\to\infty$, such that
\begin{equation*}
\lim_{k\rightarrow\infty}\sup\int_{B_{R_0}(y_{\epsilon_k})}|\tilde{w}_{a_{2k}}|^2dx=0,
\end{equation*}
where $\epsilon_k:=\epsilon(a_{1k},a_{2k})>0$ is defined by (\ref{eq2.9a}). Since $\tilde w_{a_i}$ is bounded uniformly in $H^1(\R^2)\cap L^\gamma(\R^2)$ for all $2\leq \gamma<\infty$, we may choose $\gamma>4$ and $\theta\in(0,1)$ such that $\frac{1}{4}=\frac{1-\theta}{\gamma}+\frac{\theta}{2}$. It then follows from the H\"{o}lder  inequality that
\begin{equation*}
\begin{split}
 \int_{B_{R_0}(y_{\epsilon_k})}|\tilde{w}_{a_{1k}}|^2|\tilde{w}_{a_{2k}}|^2dx&\leq
\bigg(\int_{B_{R_0}(y_{\epsilon_k})}|\tilde{w}_{a_{1k}}|^4dx\bigg)
^\frac{1}{2}\bigg(\int_{B_{R_0}(y_{\epsilon_k})}|\tilde{w}_{a_{2k}}|^4dx\bigg)^\frac{1}{2}\\
&\leq
C\bigg(\int_{B_{R_0}(y_{\epsilon_k})}|\tilde{w}_{a_{2k}}|^\gamma dx\bigg)^\frac{2(1-\theta)}{\gamma}\bigg(\int_{B_{R_0}(y_{\epsilon_k})}
|\tilde{w}_{a_{2k}}|^2dx\bigg)^\theta\\
&\leq C
\bigg(\int_{B_{R_0}(y_{\epsilon_k})}|\tilde{w}_{a_{2k}}|^2dx\bigg)^\theta\to 0
\quad \mbox{as}\quad k\to\infty .
\end{split}
\end{equation*}
By applying the estimate (\ref{eq2.14}) for $\tilde{w}_{a_{1k}}$, we use again the H\"{o}lder inequality to derive from the above that, for $k$ large,
\begin{equation*}
\begin{split}
\int_{B_{R_0}(y_{\epsilon_k})}(|\tilde{w}_{a_{1k}}|^2-|\tilde{w}_{a_{2k}}|^2)^2dx&\geq
\frac{1}{2}\int_{B_{R_0}(y_{\epsilon_k})}|\tilde{w}_{a_{1k}}|^4dx\\
&\geq
\frac{1}{2\pi
R_0^2}\bigg(\int_{B_{R_0}(y_{\epsilon_k})}|\tilde{w}_{a_{1k}}|^2dx\bigg)^2
\geq\frac{\eta_1^2}{2\pi
R_0^2},
\end{split}\end{equation*}
a contradiction, since Lemma
\ref{lem3.1} implies that
$$\inte\big(|\tilde{w}_{a_{1k}}^2-|\tilde{w}_{a_{2k}}|^2\big)^2dx=\epsilon_k^2\inte \big(|u_{a_{1k}}|^2-|u_{a_{2k}}|^2\big)^2dx\to 0
\quad \mbox{as}\quad k\to\infty .$$
Therefore, (\ref{eq2.14}) holds also for $\tilde w_{a_2}$ with some constant
$\eta_2>0$, and  part (i) is proved.

\textbf{(ii).} We first prove (\ref{eq3.20}).
By (\ref{eq2.8}) we note that
\begin{equation}\label{eq2.17}
\sum_{i=1}^2\int_{\mathbb{R}^2}V_i(x)|u_{a_i}(x)|^2dx=\sum_{i=1}^2
\int_{\mathbb{R}^2}V_i(\epsilon
x+\epsilon y_\epsilon)|w_{a_i}(x)|^2dx\to 0
\end{equation}
as $(a_1,a_2)\nearrow (a^*,a^*)$. On the contrary, suppose (\ref{eq3.20}) is  incorrect, then there exist $\delta>0$ and a subsequence $\{(a_{1n},a_{2n})\}$, which satisfies
$(a_{1n},a_{2n})\nearrow (a^*,a^*)$ as $n\to\infty $,   such that
$$\epsilon_n:=\epsilon(a_{1n},a_{2n}) \to 0 \quad \text{and}\quad {\rm dist}(\epsilon_n
y_{\epsilon_n}, \Lambda)\geq \delta>0\quad \text{as}\quad n\to\infty
.$$
This implies that there exists $C=C(\delta)>0$ such that
$$\lim_{n\to\infty}V_1(\epsilon_n
y_{\epsilon_n})+V_2(\epsilon_n
y_{\epsilon_n})\geq C(\delta)>0.$$
Hence, by  Fatou's Lemma and (\ref{eq2.16}) we obtain that
\begin{equation*}
\begin{split}
&\lim_{n\rightarrow\infty}\sum_{i=1}^2\int_{\mathbb{R}^2}V_i(\epsilon_n
x+\epsilon_n y_{\epsilon_n})|w_{a_{in}}(x)|^2dx\\
&\geq\sum_{i=1}^2
\int_{B_{R_0}(0)}\liminf_{n\rightarrow\infty}V_i(\epsilon_n
x+\epsilon_n y_{\epsilon_n})|w_{a_{in}}(x)|^2dx\geq C(\delta)\min
\{\eta_1,\eta_2\},
\end{split}
\end{equation*}
which however contradicts (\ref{eq2.17}). Therefore, (\ref{eq3.20}) holds.

%$\epsilon y_\epsilon$ is
%uniformly bounded as $\epsilon \to 0$, and moreover,  for any
%subsequence $(a_{1k},a_{2k})$  there exists a   convergent subsequence, still
%denoted by $(a_{1k},a_{2k})$, such that $\epsilon_k y_{\epsilon_k}\rightarrow
%x_0$  as $(a_{1k},a_{2k})\nearrow (a^*,a^*)$ for some point $x_0\in \R^2$. Finally,
%if $x_0\not\in \Lambda$, we then have
% $V_1(x_0)>0$ or $V_2(x_0)>0$. If $V_1(x_0)>0$ happens, then by applying Fatou's Lemma again, we deduce from (\ref{eq2.16}) that
%\begin{equation*}
%\begin{split}
%&\lim_{k\rightarrow\infty}\int_{\mathbb{R}^2}V_1(\epsilon_k
%x+\epsilon_k y_{\epsilon_k})|w_{a_{1k}}(x)|^2dx\\
%&\geq
%\int_{\mathbb{R}^2}\lim_{k\rightarrow\infty}V_1(\epsilon_k
%x+\epsilon_k y_{\epsilon_k})|w_{a_{1k}}(x)|^2dx\geq
%\frac{\eta_1}{2}V_1(x_0)
%>0,
%\end{split}\end{equation*}
%which contradicts (\ref{eq2.17}). Similarly, it is impossible to have $V_2(x_0)>0$. So we have $x_0\in\Lambda$, and Claim 2 is therefore established.\\

We prove now (\ref{eq3.21}) and (\ref{eq2.21}). Since $(u_{a_1},u_{a_2})$ is a non-negative minimizer of (\ref{eq1.3}),
it satisfies the Euler-Lagrange system
\begin{equation}\label{3:ua}
\begin{cases}
 -\Delta u_{a_1}+V_1(x)u_{a_1}=\mu_{a_1} u_{a_1}+a_1
u_{a_1}^3-\beta (u_{a_1}^2-u_{a_2}^2)u_{a_1}\quad \mbox{in}\quad \R^2\,,\\
-\Delta u_{a_2}+V_2(x)u_{a_2}=\mu_{a_2} u_{a_2}+a_2 u_{a_2}^3-\beta
(u_{a_2}^2-u_{a_1}^2)u_{a_2}\quad \mbox{in}\quad \R^2,
\end{cases}
\end{equation}
where $(\mu_{a_1},\mu_{a_2})$ is a suitable Lagrange
multiplier, and
\begin{equation*}
\mu_{a_i}=E_{a_i}^i(u_{a_i})-\frac{a_i}{2}\int_{\mathbb{R}^2}|u_{a_i}|^4dx-
\beta\int_{\R^2}(-1)^i(u_{a_1}^2-u_{a_2}^2)u_{a_i}^2dx,\quad
i=1,\,2\,.
\end{equation*}
It then follows from Lemma \ref{lem3.1}, (\ref{f}), (\ref{3:L0}) and
(\ref{eq2.10}) that for $i=1,\,2$,
\begin{equation}\label{3:mu} \mu_{a_i}\sim
-\frac{a_i}{2}\int_{\mathbb{R}^2}|u_{a_i}|^4dx\sim
-\epsilon^{-2}\quad \text{and}\quad\mu_{a_1}/\mu_{a_2}\to1\quad
\text{as}\quad (a_1,a_2)\nearrow (a^*,a^*).\end{equation}
Note also from
(\ref{eq2.15}) that $w_{a_i}(x)$ defined in (\ref{eq2.15}) satisfies the
elliptic system
\begin{equation}\label{eq2.19}
\begin{cases}
-\Delta w_{a_1}+\epsilon^2V_1(\epsilon x+\epsilon
y_{\epsilon})w_{a_1}=\epsilon^2\mu_{a_1} w_{a_1}+a_1 w_{a_1}^3-\beta
(w_{a_1}^2-w_{a_2}^2)w_{a_1}\  \mbox{ in }\ \R^2\,,\\
-\Delta w_{a_2}+\epsilon^2V_2(\epsilon x+\epsilon
y_{\epsilon})w_{a_2}=\epsilon^2\mu_{a_2} w_{a_2}+a_2
w_{a_2}^3-\beta(w_{a_2}^2-w_{a_1}^2)w_{a_2}\ \mbox{ in }\ \R^2,
\end{cases}
\end{equation}
where the Lagrange multiplier $(\mu_{a_1},\mu_{a_2})$ satisfies (\ref{3:mu}).

For any given  sequence $\{(a_{1k},a_{2k})\}$
with $(a_{1k},a_{2k})\nearrow (a^*,a^*)$ as $k\to\infty$, we deduce  from (\ref{eq3.20}) and (\ref{3:mu}) that there exists a  subsequence of $\{(a_{1k},a_{2k})\}$,
still denoted by $\{(a_{1k},a_{2k})\}$, such that
\begin{equation}\label{eq3.28}
\epsilon_ky_{\epsilon_k} \overset{k}{\to} x_0\in \Lambda,\ \ \epsilon_k^2\mu_{a_{ik}}\overset{k}{\to}-\lambda^2 < 0\ \, \text{  for some } \ \lambda>0,\end{equation}
 and
$$w_{a_{ik}}\overset{k}{\rightharpoonup}
w_i\geq 0\text{ weakly  in } H^1(\mathbb{R}^2)\text{ for some } w_i\in
H^1(\R^2), \ \, i=1,\,2.$$
Since (\ref{3:est}) implies that
\begin{equation}\label{3:two}\|w_{a_1}^2-w_{a_2}^2\|_2=\epsilon\|u_{a_1}^2-u_{a_2}^2\|_2\to0 \quad {\rm as}\quad (a_1,a_2)\nearrow
(a^*,a^*),\end{equation}
we have $w_1=w_2\geq0$ a.e. in $\R^2$. We thus  write $0\leq w_0:=w_1=w_2\in H^1(\R^2)$. Passing to the weak limit
in (\ref{eq2.19}), we deduce from (\ref{eq3.28}) and (\ref{3:two}) that $w_0$ satisfies
\begin{equation}\label{eq2.20}
-\Delta w_0(x)=-\lambda^2 w_0(x)+a^* w_0^{3}(x)\ \, \text {in}
\ \, \mathbb{R}^2.
\end{equation}
Furthermore, it follows from (\ref{eq2.16}) and the strong maximum principle that $w_0>0$. By a
simple rescaling,  the uniqueness (up to translations) of positive
solutions for the equation (\ref{Kwong})
implies that
\begin{equation}\label{3.34'}
w_0(x)=\frac{\lambda}{\|Q\|_2}Q(\lambda|x-y_0|) \ \, \text{for
some} \ \, y_0\in\mathbb{R}^2.
\end{equation}
Note that $||w_0||_2^2=1$, by the norm preservation we further conclude
that $$w_{a_{ik}}\overset{k}{\to} w_0\ \text{ strongly in $L^2(\mathbb{R}^2)$}.$$
Moreover, this strong convergence  holds also for all $p\geq2$, since $\{w_{a_{ik}}\}$ is bounded in $H^1(\R^2)$. Then, note that $w_{a_{ik}}$ and $w_0$ satisfy (\ref{eq2.19}) and (\ref{eq2.20}),
respectively,   a simple analysis shows that
$$w_{a_{ik}}\overset{k}{\to} w_0\ \text{ strongly in $H^1(\mathbb{R}^2)$},\ \ i=1,\,2.$$
Therefore, (\ref{eq3.21}) and (\ref{eq2.21}) are established.\qed

Applying above lemmas, we finally  prove Proposition \ref{3:thm2} on the optimal estimates of $e(a_1,a_2)$.

\noindent{\bf Proof of Proposition  \ref{3:thm2}.} For any sequence $\{(a_{1k},a_{2k})\}$
satisfying $(a_{1k},a_{2k})\nearrow (a^*,a^*)$ as $k\to\infty$,  it follows from Lemma \ref{le3.3} (ii) that there exists a convergent subsequence,
still denoted by $\{(a_{1k},a_{2k})\}$, such that (\ref{eq3.21}) holds and  $\epsilon_k
y_{\epsilon_k} \overset{k}\to x_0\in \Lambda$. Without loss of generality, we may assume $x_0=x_{1j_0}$ for some $1\leq j_0\leq l$. We first claim that
\begin{equation}\label{eq3.31}\limsup_{k\to\infty}\frac{|\epsilon_k
y_{\epsilon_k}-x_{1j_0}|}{\epsilon_k}<\infty.\end{equation}
Actually, by (\ref{eq3.21}) and (\ref{eq2.21}), we have for some $R_0>0$,
 \begin{equation}\label{3:lat}
 \begin{split}
 e(a_{1k},a_{2k})
  =&E_{a_{1k},a_{2k}}(u_{a_{1k}},u_{a_{2k}})\\
 \geq&\sum_{i=1}^2\Big\{\frac{1}{\epsilon
_k^2}\Big[\int_{\mathbb{R}^2}|\nabla
w_{a_{ik}}(x)|^2dx-\frac{a^*}{2}\int_{\mathbb{R}^2}|w_{a_{ik}}(x)|^4dx\Big]\\
&+\frac{a^*-{a_{ik}}}{2\epsilon
_k^2}\int_{\mathbb{R}^2}|w_{a_{ik}}(x)|^4dx+
\int_{\mathbb{R}^2}V_i(\epsilon _k
x+\epsilon _ky_{\epsilon_k})|w_{a_{ik}}(x)|^2dx\Big\}\\
 \geq &\sum_{i=1}^2\Big\{\frac{a^*-{a_{ik}}}{4\epsilon
_k^2}\int_{\mathbb{R}^2}|w_0(x)|^4dx+
\int_{B_{R_0}(0)}V_i(\epsilon _k
x+\epsilon _ky_{\epsilon_k})|w_{a_{ik}}(x)|^2dx\Big\}\\
 \geq &
C_1\frac{2a^*-a_{1k}-a_{2k}}{\epsilon_k^2}+C_2\sum_{i=1}^2
\epsilon_k^{p_{ij_0}}\int_{B_{R_0}(0)}\Big|x+
\frac{\epsilon_ky_{\epsilon_k}-x_{1j_0}}{\epsilon_k}\Big|^{p_{ij_0}}|w_{a_{ik}}(x)|^2dx.
\end{split}
\end{equation}
Suppose now that there exists a subsequence such that   $\frac{|\epsilon_k
y_{\epsilon_k}-x_{1j_0}|}{\epsilon_k}\to\infty$ as $k\to\infty$. By  using Fatou's Lemma, we then deduce from (\ref{3:lat}) and (\ref{eq2.16}) that for any $M>0$,
\begin{equation*}
\begin{split}
e(a_{1k},a_{2k})&\geq
C_1\frac{2a^*-a_{1k}-a_{2k}}{\epsilon_k^2}+C_2M\epsilon_k^{\bar
p_{j_0}}\geq CM^{\frac{1}{\bar
p_{j_0}+2}}\Big(a^*-\frac{a_{1k}+a_{2k}}{2}\Big)^\frac{\bar
p_{j_0}}{\bar p_{j_0}+2}\\
&\geq CM^{\frac{1}{\bar
p_{j_0}+2}}\Big(a^*-\frac{a_{1k}+a_{2k}}{2}\Big)^\frac{ p_0}{
p_0+2},
\end{split}
\end{equation*}
where $p_0=\max_{1\leq j\leq
l}\bar p_j$ and $\bar
p_{j_0}=\min\{p_{1j_0},p_{2j_0}\}>0$    are given by (\ref{1:def}).
This estimate however contradicts the upper bound of (\ref{3:est}), and the claim (\ref{eq3.31}) is thus established.

We now deduce from (\ref{eq3.31}) that there exists a subsequence of $\{\epsilon_k\}$, still denoted by $\{\epsilon_k\}$,
such that
$$\frac{\epsilon_ky_{\epsilon_k}-x_{1j_0}}{\epsilon_k}\to y_1 \ \,\text{as}\ \, k\to\infty$$
holds for some $y_1 \in \R^2$. By applying (\ref{3:lat}), then there exists a constant $C_1>0$, independent of $a_{1k}$ and $a_{2k}$, such that
$$e(a_{1k},a_{2k})\geq
C_1\Big(a^*-\frac{a_{1k}+a_{2k}}{2}\Big)^\frac{\bar p_{j_0}}{\bar
p_{j_0}+2}  \ \  \text{as}\ \,
(a_{1k},a_{2k})\nearrow (a^*,a^*).$$
Since $\bar p_{j_0}\leq p_0$, applying the upper bound of (\ref{3:est}), we conclude from the above estimate that $\bar p_{j_0}= p_0$ and (\ref{energy}) holds for the above subsequence $\{(a_{1k},a_{2k})\}$.

We next prove that (\ref{1:u4}) holds for the above subsequence
$\{(a_{1k},a_{2k})\}$. Suppose that
$$ \mbox{either} \quad \epsilon_k>>
\Big(a^*-\frac{a_{1k}+a_{2k}}{2}\Big)^\frac{1}{p_0+2}\quad  \mbox{or}
\quad 0<\epsilon_k<<
\Big(a^*-\frac{a_{1k}+a_{2k}}{2}\Big)^\frac{1}{p_0+2}\quad \mbox{as}\quad k\to\infty,$$
it then follows from (\ref{3:lat}) that $e(a_{1k},a_{2k})>>
\Big(a^*-\frac{a_{1k}+a_{2k}}{2}\Big)^\frac{ p_0}{p_0+2}$ as $k\to\infty$, which however contradicts (\ref{energy}). This completes the proof of (\ref{1:u4}).

Since the above argument shows that Proposition \ref{3:thm2} holds for any given subsequence
$\{(a_{1k},a_{2k})\}$ with $(a_{1k},a_{2k})\nearrow (a^*,a^*)$, an approach similar to that of \cite{GZZ} then gives that Proposition
\ref{3:thm2} holds essentially for the whole sequence $\{(a_1,a_2)\}$ satisfying $(a_1,a_2)\nearrow (a^*,a^*)$.  \qed

\section{Proof of Theorem \ref{thm1.4}}

This section is devoted to the proof of Theorem \ref{thm1.4} on the mass concentration of minimizers.
Using the  notations as in (\ref{f})-(\ref{eq1.3}), in order to prove Theorem \ref{thm1.4} on the minimizers of (\ref{eq1.4}) as $(b_1+\beta,b_2+\beta)\nearrow(a^*,a^*)$, it suffices to establish the following theorem on the minimizers of (\ref{eq1.3}) as $(a_1,a_2)\nearrow(a^*,a^*)$.
\begin{thm}\label{4:thm3}
Assume that $0<\beta<a^*$ and  $V_i(x)$ satisfies (\ref{1:V}) and (\ref{1:V1}) for $i=1$ and $2$.
For any  sequence
$\{(a_{1k},a_{2k})\}$ satisfying $(a_{1k},a_{2k})\nearrow (a^*,a^*)$ as
$k\to\infty$, and let $(u_{a_1k},u_{a_2k})$ be the corresponding non-negative
minimizer of (\ref{eq1.3}).
Then there exists a subsequence of $\{(a_{1k},a_{2k})\}$, still denoted by
$\{(a_{1k},a_{2k})\}$, such that  each $u_{a_{ik}}$ ($i=1$, $2$) has a unique global maximum point $x_{ik}$ satisfying
\begin{equation}\label{1:rate1}
x_{ik} \longrightarrow \bar x_0\in \mathcal{Z} \text{ and }
\frac{|x_{ik}-\bar x_0|}{\big(a^*-\frac{a_{1k}+a_{2k}}{2}\big)^\frac{1}{p_0+2}}\longrightarrow 0 \quad
\text{as}\ \, k\to\infty.
\end{equation}
Moreover, for $i=1$ and $2$,
\begin{equation*}
\lim_{k\to\infty} \Big(a^*-\frac{a_{1k}+a_{2k}}{2}\Big)^\frac{1}{p_0+2} u_{a_{ik}}\Big(\big(a^*-\frac{a_{1k}+a_{2k}}{2}\big)^\frac{1}{p_0+2}
x+x_{ik}\Big)=\frac{\lambda}{\|Q\|_{2}}Q(\lambda x)
 \end{equation*}
 strongly in $H^1(\R^2)$, where $\lam >0$ is given by
\begin{equation}\label{def:lam}
  \lambda = \Big(
\frac{p_0\gamma}{4}\inte | x|^{p_0}Q^2(x)dx\Big)^\frac{1}{p_0+2},
\end{equation}
$p_0>0$ and  $\gamma >0$  are defined in (\ref{1:def})-(\ref{def:Z}).
\end{thm}
 Let  $(u_{a_1},u_{a_2})$  be a
non-negative minimizer of (\ref{eq1.3}), where $(a_1,a_2)\nearrow (a^*,a^*)$. Define
\begin{equation}\label{3:eps}
 \eps :=\big(a^*-\frac{a_1+a_2}{2}\big)^{\frac{1}{p_0+2}}>0.
\end{equation}
It then follows from (\ref{1:u4}) and (\ref{3:est}) that, as $(a_1,a_2)\nearrow (a^*,a^*)$,
\begin{equation}\label{4:V}
\sum_{i=1}^2\int_{\mathbb{R}^2}V_i(x)|u_{a_i}(x)|^2dx\leq
e(a_1,a_2)<C\big(a^*-\frac{a_1+a_2}{2}\big)^\frac{p_0}{p_0+2}
\end{equation}
and
\begin{equation}\label{4:ua}
\int_{\mathbb{R}^2}|\nabla u_{a_i}(x)|^2dx\sim \eps
^{-2},\quad\int_{\mathbb{R}^2}| u_{a_i}(x)|^4dx\sim \eps ^{-2},
\end{equation}
where $i=1,\,2$. Similar to (\ref{eq2.16}), we know  that there
exist a sequence $\{y_\varepsilon\}$ as well as positive constants $R_0$
and $\eta_i$  such that
\begin{equation}\label{3:eq2.16}
\liminf_{\varepsilon\searrow
0}\int_{B_{R_0}(0)}|w_{a_i}|^2dx\geq\eta_i>0,\quad i=1,\,2,
\end{equation}
where $w_{a_i}$ is the $L^2(\mathbb{R}^2)$-normalized
function defined by \begin{equation}\label{4:wa}
 w_{a_i}(x)=\varepsilon
u_{a_i}(\varepsilon x+\varepsilon y_\varepsilon),\quad i=1,\,2.
\end{equation}
Note from (\ref{4:ua}) that
 \begin{equation}\label{3:4}
  M\leq\int_{\mathbb{R}^2}|\nabla w_{a_i}|^2dx\leq \frac{1}{M},\quad
M\leq\int_{\mathbb{R}^2}| w_{a_i}|^4dx\leq \frac{1}{M},\quad i=1,\,2,
\end{equation}
where the positive constant $M$ is independent of $a_1$ and $a_2$.\\

\noindent\textbf{Proof of Theorem \ref{4:thm3}.} Let $\eps_k
:=\big(a^*-\frac{a_{1k}+a_{2k}}{2}\big)^{\frac{1}{p_0+2}}>0$, where
$(a_{1k},a_{2k})\nearrow (a^*,a^*)$ as $k\to\infty$, and denote
$(u_{1k}(x),u_{2k}(x)):=(u_{a_{1k}}(x),u_{a_{2k}}(x))$ a
non-negative minimizer of (\ref{eq1.3}). Inspired by
\cite{GS,Wang}, we shall complete the proof of Theorem
\ref{4:thm3} by the following three steps.
\\

\noindent{\em Step 1: Decay estimate for $(u_{1k}(x),u_{2k}(x))$.} Let
$w_{ik}(x):=w_{a_{ik}}(x)\ge 0$ be defined by (\ref{4:wa}). By a similar
analysis to the proof of Lemma \ref{le3.3}(ii), we know  that there exists a
subsequence of $\{\varepsilon_k\}$, still denoted by $\{\varepsilon_k\}$,  such that
\begin{equation}\label{4.7}
z_k:=\varepsilon_ky_{\varepsilon_k}\xrightarrow{k}x_0\ \text{ for some }\ x_0\in\Lambda,\end{equation}
where the set $\Lambda$ is defined by (\ref{1:com}).
and  $w_{ik} $ $(i=1,2)$ satisfies
\begin{equation}\label{4:wa1}
\begin{cases}
-\Delta w_{1k}+\varepsilon_k^2V_1(\varepsilon_k x+z_k)w_{1k}=\varepsilon_k^2\mu_{1k} w_{1k}+a_{1k}
w_{1k}^3-\beta(w_{1k}^2-w_{2k}^2)w_{1k}\ \mbox{ in }\ \R^2\,,\\[2mm]
-\Delta w_{2k}+\varepsilon_k^2V_2(\varepsilon_k x+z_k)w_{2k}=\varepsilon_k^2\mu_{2k} w_{2k}+a_{2k}
w_{2k}^3-\beta (w_{2k}^2-w_{1k}^2)w_{2k}\ \mbox{ in }\ \R^2,
\end{cases}
\end{equation}
where $(\mu_{1k},\mu_{2k})$ is a suitable Lagrange multiplier satisfying
$$\mu_{ik}\sim -\varepsilon_k^{-2}\ \, \mbox{and}\ \, \mu_{1k}/\mu_{2k}\to 1 \quad \text{as}\quad (a_{1k},a_{2k})\nearrow (a^*,a^*),\quad i=1,\,2.$$
Moreover, $w_{ik}\xrightarrow{k} w_0$ strongly in $H^1(\mathbb{R}^2)$
for some  $w_0>0$ satisfying
\begin{equation}\label{3:eq2.20}
-\Delta w_0(x)=-\lambda ^2 w_0(x)+a^* w_0^{3}(x)\quad \text { \ in }
\ \mathbb{R}^2,
\end{equation}
where $\lambda >0$ is a  constant. Similar to (\ref{3.34'}), we know  that
\begin{equation}\label{4:w0}
w_0(x)=\frac{\lambda}{\|Q\|_2}Q(\lambda|x-y_0|) \ \, \text{for
some} \ \, y_0\in\mathbb{R}^2.
\end{equation}
By the exponential decay of $Q$, then  for any $\alpha>2$,
\begin{equation}\label{unif}
\int_{|x|\geq R}|w_{ik}|^\alpha dx\rightarrow0 \ \, \text{as}\ \,
R\rightarrow\infty \ \, \text{uniformly for large }  k,
\ \, \text{where}\ \, i=1,\,2.
\end{equation}
Recall from (\ref{4:wa1}) that $-\Delta w_{ik}(x)\leq c_i(x)w_{ik}(x)$ in $\R^2$,
where $c_i(x)=a_{ik} w_{ik}^2+(-1)^i\beta(w_{1k}^2-w_{2k}^2)$ in $\R^2$ for $i=1,\, 2$. By
applying De Giorgi-Nash-Moser theory, we thus have
$$\max_{B_1(\xi)} w_{ik}\leq C\Big(\int_{B_2(\xi)}|w_{ik}|^\alpha
dx\Big)^\frac{1}{\alpha}, \ \ i=1,\, 2,$$
where $\xi$ is an arbitrary point in
$\mathbb{R}^2$, and $C>0$ is a constant depending only on the bound of
$\|w_{1k}\|_{L^\alpha(B_2(\xi))}+\|w_{2k}\|_{L^\alpha(B_2(\xi))}$.
Hence  (\ref{unif}) implies  that
\begin{equation}
w_{ik}(x)\rightarrow 0 \quad \text{as}\quad |x|\rightarrow\infty
\quad \text{uniformly for } k,\quad i=1,\,2.
\end{equation}
Since $w_{ik}$ ($i=1,\,2$) satisfies (\ref{4:wa1}), apply the comparison
principle as in \cite{KW} to compare $w_{ik}$ with
$Ce^{-\frac{\lambda}{2}|x|}$ ($\lambda >0$ is obtained in (\ref{3:eq2.20})),  which then shows that there exist a constant $C>0$ and a
large constant $R>0$, independent of $k$, such that
\begin{equation}\label{exp}
w_{ik}(x)\leq Ce^{-\frac{\lambda}{2}|x|} \quad \text{for} \quad
|x|>R \quad \text{as}\quad k\rightarrow\infty,\quad i=1,\,2.
\end{equation}
\vskip 0.1truein

\noindent{\em Step 2: The detailed concentration behavior.} For the
convergent subsequence  $\{w_{ik}(x)\}$ obtained in
Step 1, let
$\bar{z}_{ik}$ be any local maximum point of $u_{ik}(x),$ ($i=1,\,2$). We
claim that there is a sequence  $\{k\}$, passing to a subsequence if necessary,  such that
\begin{equation}\label{4:zi}\lim_{k\to\infty}\frac{|\bar z_{1k}-\bar
z_{2k}|}{\varepsilon_k}=0\,.\end{equation}

For showing (\ref{4:zi}), we first prove that
\begin{equation}\label{4:z1}\limsup_{k\to\infty}\frac{|\bar
z_{ik}-z_k|}{\varepsilon_k}<\infty,\quad i=1,\,2,\end{equation}
where
$z_k=\varepsilon_ky_{\varepsilon_k} x_0\in\Lambda$  as $k\to\infty$.
Indeed, if (\ref{4:z1}) is false, i.e.,  $\frac{|\bar
z_{ik}-z_k|}{\varepsilon_k} \stackrel{k} \rightarrow \infty$  holds for $i=1$ or
$2$. Without loss of generality, we assume that $\frac{|\bar
z_{1k}-z_k|}{\varepsilon_k} \stackrel{k} \rightarrow \infty$. It follows from (\ref{4:wa}) and (\ref{exp})  that
$$u_{ik}(\bar z_{1k})=\frac{1}{\eps_k}w_{ik}\big(\frac{\bar z_{1k}-z_k}{\eps_k}\big)=o\big(\frac{1}{\eps_k}\big)\quad \mbox{as}\quad k\to\infty,\quad i=1,\,2.$$
This however leads to a contradiction, since  (\ref{3:ua}) implies
that
$$a_{1k}u_{1k}^2(\bar z_{1k})-\beta\big(u_{1k}^2(\bar z_{1k})-u_{2k}^2(\bar z_{1k})\big)\geq -\mu_{1k}\geq C\varepsilon_k^{-2}.$$
%$i.e.$, $\varepsilon_ku_{ik}(\bar z_{1k})>C$ uniformly for some $C>0$.
Therefore, (\ref{4:z1}) is true.

By (\ref{4:z1}),  there exists a
sequence  $\{k\}$ such that
$$\lim_{k\to\infty}\frac{\bar
z_{ik}-z_k}{\varepsilon_k}=y_i\quad \text{for some}\ \, y_i\in
\R^2\,,\ \, i=1,\,2.$$
Set
\begin{equation}\label{4:wb}
\bar
w_{ik}(x):=w_{ik}\big(x+\frac{\bar
z_{ik}-z_k}{\varepsilon_k}\big)=\eps_ku_{ik}\big(\eps_k x+\bar z_{ik}\big),\quad i=1,\,2.
\end{equation}
By  Step 1, $w_{ik}\to w_0$ strongly in $H^1(\R^2)$ as $k\to\infty$, and $w_0>0$ satisfying (\ref{4:w0}),
then
\begin{equation*}\lim_{k\to\infty}\bar w_{ik}(x)=w_0(x+y_i)=\frac{\lambda}{\|Q\|_2}Q(\lambda|x+y_i-y_0|)\ \, \text{strongly in }H^1(\R^2), \ \, i=1,\,2.\end{equation*}
Since  the origin $(0,0)$ is a critical point of $\bar{w}_{ik}$  for
all $k>0$ ($i=1,\,2$), it is also a critical point of $w_0(x+y_i)$. On the other hand, $Q(\lambda|x-z_0|)$ possesses $z_0$ as its unique critical (maximum) point. We therefore conclude  that   $w_0(x+y_i)=\frac{\lambda}{\|Q\|_2}Q(\lambda|x+y_i-y_0|)$ is
spherically symmetric about the origin. Hence, $y_i=y_0$,
and
\begin{equation}\label{eq2.350}
\lim_{k\to\infty}\bar
w_{ik}(x)=\frac{\lambda}{\|Q\|_2}Q(\lambda|x|):=\bar w_0 \ \  \text{strongly in }H^1(\R^2),\ \, i=1,\,2.
\end{equation}
The estimate (\ref{4:zi}) is followed by   (\ref{4:wb}) and (\ref{eq2.350}).

Similar to the discussion of  proving (\ref{4:z1}), we know  that  each
local maximum point of $\bar{w}_{ik}(x)$ ($i=1,\,2$), which is also a local maximum
point of $u_{ik}(x)$, must stay in a finite ball in $\mathbb{R}^2$.  Since $V_i(x)\in C_{\rm loc}^\alpha(\R^2)$, we can deduce from (\ref{4:wa1}) and  standard elliptic regular theory   that
$${w}_{ik}(x)\xrightarrow{k} {w}_0(x)\ \text{ in }\
C^{2,\alpha}_{loc}(\mathbb{R}^2),\ \, i=1,\,2.$$
Moreover, by (\ref{4:wb}) and  (\ref{4:zi}), we can further obtain that
$$\bar{w}_{ik}(x)\xrightarrow{k} \bar{w}_0(x)\ \text{ in } \
C^{2,\alpha}_{loc}(\mathbb{R}^2),\ \, i=1,\,2.$$
  Because
the origin $(0,0)$ is the only non-degenerate  critical point of
$\bar{w}_0(x)$, all local maximum points of $\bar{w}_{ik}(x)$ must approach the origin and
 stay in a small ball $B_\epsilon(0)$ as $k\rightarrow\infty$, where $\epsilon >0$ is small.
It then follows from Lemma 4.2 in \cite{NT} that for large $k$,
$\bar{w}_{ik}(x)$ has no critical points other than the origin. This gives
the uniqueness of local maximum points
for $\bar{w}_{ik}(x)$ and $u_{ik}(x)$ ($i=1,\,2$).\\

\noindent{\em Step 3: Completion of the proof.}
Let
\begin{equation}\label{4:def}
 \gamma_j(x) =  \frac{V_1(x)+V_2(x)}{|x-x_{1j}|^{\bar p_j}}, \quad 1\leq j\leq l,
\end{equation}
where $\bar p_j>0$ is defined by (\ref{1:def}), so that the limit
$\lim_{x\to x_{1j}}\gamma_j(x)=\gamma_j(x_{1j})$  exists
for all $i\in\{1,\cdots ,l\}$ in view of the assumptions on $V_1$ and $V_2$. Moreover, one can note that
$\gamma_j(x_{1j})=\gamma_j\geq \gamma$ if $x_{1j}\in \bar \Lambda$,
where $\gamma_j$ and $\bar\Lambda$ are defined by (\ref{def:li}) and
(\ref{1:ba}), respectively. We hence obtain  from (\ref{4:wb}) that
 \begin{equation}\label{4:lat}
 \begin{split}
e(a_{1k},a_{2k})& =E_{a_{1k},a_{2k}}(u_{a_{1k}},u_{a_{2k}}) \\
&\geq\sum_{i=1}^2\Big\{\frac{1}{\eps
_k^2}\Big[\int_{\mathbb{R}^2}|\nabla
\bar w_{{ik}}(x)|^2dx-\frac{a^*}{2}\int_{\mathbb{R}^2}|\bar w_{{ik}}(x)|^4dx\Big] \\
&+\frac{a^*-{a_{ik}}}{2\eps _k^2}\int_{\mathbb{R}^2}|\bar
w_{{ik}}(x)|^4dx+ \int_{\mathbb{R}^2}V_i(\eps _k
x+\bar z_{ik})|\bar w_{{ik}}(x)|^2dx\Big\} \\
&\geq\sum_{i=1}^2\bigg[\frac{a^*-{a_{ik}}}{2\eps
_k^2}\int_{\mathbb{R}^2}|\bar w_{{ik}}(x)|^4dx+
\int_{\mathbb{R}^2}V_i(\eps_k
x+\bar z_{ik})|\bar w_{{ik}}(x)|^2dx\bigg],
\end{split}\end{equation}
where $\bar{z}_{ik}$ is the unique global maximum point of $u_{ik}$,
and $\bar{z}_{ik} \overset{k}\to x_0\in \Lambda$, $i=1,\,2$. We may assume
that $x_0=x_{1j_0}$ for some $1\leq j_0\leq l$.

We first  claim that
\begin{equation}\label{4:cla} \frac{|\bar{z}_{ik}-x_{1j_0}|}{\varepsilon_k}
\quad \text{is uniformly bounded  as}\ \, k\to\infty,\quad
\text{where}\ \, i=1,\,2.\end{equation}
Otherwise, if (\ref{4:cla}) is false for $i=1$ or $2$. It then follows from (\ref{4:zi}) that both of them
are unbounded, and hence there exists a subsequence of $\{(a_{1k},a_{2k})\}$, still
denoted by $\{(a_{1k},a_{2k})\}$,  such that
$$\lim_{k\to\infty}\frac{|\bar{z}_{ik}-x_{1j_0}|}{\varepsilon_k}=\infty,\quad  i=1,\,2.$$
We then derive from Fatou's Lemma that for any $M>0$ large enough,
 \begin{equation*}
\begin{split}
 &\lim_{k \to \infty}\inf \eps_k^{-\bar p_{j_0}}\sum_{i=1}^2
\int_{\mathbb{R}^2}V_i(\eps _k x+\bar z_{ik})|\bar
w_{{ik}}(x)|^2dx \\
=&\lim_{k \to \infty} \inf\sum_{i=1}^2
\int_{\mathbb{R}^2}\frac{V_i(\eps _k x+\bar
z_{ik})}{|\eps_kx+\bar z_{ik}-x_{1j_0}|^{\bar
p_{j_0}}}\Big|x+\frac{\bar
z_{ik}-x_{1j_0}}{\eps_k}\Big|^{\bar p_{j_0}}|\bar w_{{ik}}(x)|^2dx\\
\geq &\sum_{i=1}^2
\int_{\mathbb{R}^2}\lim_{k \to \infty}\inf\Big(\frac{V_i(\eps _k x+\bar
z_{ik})}{|\eps_kx+\bar z_{ik}-x_{1j_0}|^{\bar
p_{j_0}}}\Big|x+\frac{\bar
z_{ik}-x_{1j_0}}{\eps_k}\Big|^{\bar p_{j_0}}|\bar w_{{ik}}(x)|^2\Big)dx\geq M.
\end{split}
\end{equation*}
This estimate and (\ref{4:lat})  imply that
\begin{equation}\label{4:latt}e(a_{1k},a_{2k})\geq M \eps_k^{\bar p_{j_0}}=M
\Big(a^*-\frac{a_{1k}+a_{2k}}{2}\Big)^{\bar p_{j_0}} \end{equation}
holds for arbitrary constant $M>0$, which however contradicts
Proposition \ref{3:thm2}, due to the fact that $\bar p_{j_0}\le p_0$. Therefore, (\ref{4:cla}) is proved.

We next show that $\bar p_{j_0}=p_0$, i.e., $x_{1j_0}\in\bar\Lambda$, where the set $\bar\Lambda$ is defined by (\ref{1:ba}).
By  (\ref{4:cla}), we know  that  there exists a
subsequence of $\{(a_{1k},a_{2k})\}$ such that
\begin{equation}\label{eq2.46a}
\lim_{k\to\infty}\frac{\bar{z}_{ik}-x_{1j_0}}{\varepsilon_k}=\bar
z_0\quad\text{for some }\ \bar{z}_0\in\mathbb{R}^2,\quad i=1,\,2.
\end{equation}
Since $Q$ is a radially
decreasing function and decays exponentially as $|x|\to\infty$, we
then deduce from (\ref{eq2.350}) that
 \begin{equation}
 \begin{split}
 &\lim_{k \to \infty}\inf \eps_k^{-\bar p_{j_0}}\sum_{i=1}^2
\int_{\mathbb{R}^2}V_i(\epsilon _k x+\bar z_{ik})|\bar
w_{{ik}}(x)|^2dx \\
=&\lim_{k \to \infty}\inf \sum_{i=1}^2
\int_{\mathbb{R}^2}\frac{V_i(\epsilon _k x+\bar
z_{ik})}{|\eps_kx+\bar z_{ik}-x_{1j_0}|^{\bar
p_{j_0}}}\Big|x+\frac{\bar
z_{ik}-x_{1j_0}}{\eps_k}\Big|^{\bar p_{j_0}}|\bar w_{{ik}}(x)|^2dx \\
\ge &\,\gamma_{j_0}({x_{1j_0}})\inte|x+\bar z_0|^{\bar p_{j_0}}\bar
w_0^2dx
 \geq
 \frac{\gamma_{j_0}({x_{1j_0}})}{\lambda^{\bar p_{j_0}}\|Q\|_2^2}\inte|x|^{\bar p_{j_0}}Q^2dx,\end{split}\label{4:lat1}
\end{equation}
where the equality holds if and only if $\bar z_0=(0,0)$. We hence obtain from (\ref{4:lat}) and (\ref{4:lat1}) that $\bar p_{j_0}=p_0$, otherwise we get that (\ref{4:latt}) holds with $M$ being replaced by  some
$C>0$, which however contradicts Proposition \ref{3:thm2}.

By the fact $\bar p_{j_0}=p_0$, we now have $x_{1j_0}\in \bar \Lambda$ and
$\gamma_{j_0}(x_{1j_0})=\gamma_{j_0}$. It then follows from (\ref{3:eps}),
(\ref{4:lat}) and (\ref{4:lat1}) as well as  (\ref{1:id})that
\begin{align}
 \lim_{k\to \infty}\inf\frac{e(a_{1k},a_{2k})}{\eps_k^{ p_0}} &\geq
 \|\bar w_0\|_4^4 + \gamma_{j_0}\inte|x+\bar z_0|^{ p_0}\bar w_0^2dx  \nonumber\\
&\geq\frac{1}{a^*} \Big( 2\lambda^2 + \frac{\gamma_{j_0}}{\lambda^{
p_0}}\inte|x|^{ p_0}Q^2dx\Big),\label{mb}
\end{align}
and``=" holds in the last inequality if and only if  $\bar z_0=(0,0)$.
%Moreover, since $(u_{1k}(x),u_{2k}(x))$ is a minimizer of $e(a_{1k},a_{2k})$, (\ref{4:lat1}) implies that the equality of (\ref{mb}) holds and $\bar z_0=(0,0)$.
The estimate (\ref{1:rate1}) then follows from (\ref{eq2.46a}) and this conclusion.
Further, taking the
infimum of (\ref{mb}) over $\lambda>0$ yields that
\begin{equation}\label{lim}
\lim_{k\to
\infty}\inf\frac{e(a_{1k},a_{2k})}{\eps_k^{
p_0}}\geq\frac{ 2p_0+4 }{ p_0a^*}\Big(
\frac{p_0\gamma_{j_0}\inte|x|^{ p_0}Q^2dx}{4}\Big)^\frac{2}{
p_0+2},
\end{equation}
where the equality is achieved at  $$\lam= \lam _0:=\Big(
\frac{p_0\gamma_{j_0}\inte|x|^{ p_0}Q^2dx}{4}\Big)^\frac{1}{
p_0+2}.$$

We finally remark that the limit of (\ref{lim}) actually exists and is equal to the right hand  side of (\ref{lim}). To see  this fact, we simply take
$$u_1(x)=u_2 (x)=\frac{\alpha }{\eps_k \|Q\|_2}Q\Big(\frac{\alpha
|x-x_{1j_1}|}{\eps_k}\Big)\quad \text{with}\quad
x_{1j_1}\in\mathcal{Z}$$
as a trial function for $E_{a_{1k},a_{2k}}(\cdot,\cdot)$
and minimizes it over $\alpha
>0$, which then leads to the limit
\begin{equation}\label{4:lat3}\lim_{k\to
\infty}\inf\frac{e(a_{1k},a_{2k})}{\eps_k^{
p_0}}\leq\frac{ 2p_0+4 }{ p_0a^*}\Big( \frac{p_0\gamma\inte|x|^{
p_0}Q^2dx}{4}\Big)^\frac{2}{ p_0+2}.
\end{equation}
Since $\gamma = \min\big\{\gamma_1,\ldots
,\gamma_l\big\}$, it follows from (\ref{lim}) and (\ref{4:lat3}) that  $\gamma=
\gamma_{j_0}$, i.e., $x_{1j_0}\in\mathcal{Z}$, and further (\ref{lim})
is indeed an equality. This yields that
$\lambda$ is unique, which is independent of the choice of the
subsequence, and minimizes (\ref{mb}),
i.e., $\lam=\lam _0$. Moreover, when (\ref{lim}) becomes an equality, which implies that (\ref{mb}) is also an equality. Thus, $\bar z_0=(0,0)$, which together with (\ref{eq2.46a}) give (\ref{1:rate1}).
 \qed\\

%\vspace {.95cm}

\appendix
%\appendixpage
%\addappheadtotoc

\section{Appendix: Some Proofs}

In this appendix we shall establish the following lemma and provide a different proof of Theorem \ref{thm1}.

\begin{lem}\label{lemA}
Suppose that positive constants $b_1$, $b_2$ and $\beta$ satisfy $0<b_1<b_2<a^*$ and $0<b_2\leq 2\beta+b_1$. Define
$$l(t)=\frac{t^2+t}{\frac{b_2}{2}t^2+\beta t+\frac{b_1}{2}}, \  \, t\in[0,\infty).$$
 Then we have $$l(t)> l(1)=\frac{2}{\frac{b_2}{2}+\beta+\frac{b_1}{2}},\, \ t\in (1,\infty).$$
\end{lem}

\noindent\textbf{Proof.} Direct calculation shows that
\begin{equation}
l'(t)=\frac{(\beta-\frac{b_2}{2})t^2+b_1t+\frac{b_1}{2}}{\big(\frac{b_2}{2}t^2+\beta t+\frac{b_1}{2}\big)^2}.
\end{equation}
Let
$$m(t)=(\beta-\frac{b_2}{2})t^2+b_1t+\frac{b_1}{2}, \ t\in[0,\infty).$$
If $\beta\geq \frac{b_2}{2}$, then $m'(t)=2(\beta-\frac{b_2}{2})t+b_1>0$. This implies that $l'(t)>0$ for $t\in[0,\infty)$, and we are done.

We now suppose that $\beta< \frac{b_2}{2}$. In this case, we have
$$m(t)>0  \ \, \text{ if }\ \, t\in(0, \hat t\,); \ \, m(t)<0
\ \, \text{ if }\ \, t\in(\hat t, \infty),$$
where
$$\hat t=\frac{b_1+\sqrt{b_1(b_1+b_2-2\beta)}}{b_2-2\beta}>1, $$
in view of the assumption that $b_2\leq 2\beta+b_1$.
Thus, $l(t)$ is strictly increasing in $(1,\hat t\ ]$ and strictly decreasing in $(\hat t, \infty)$.
Since $b_2\leq 2\beta+b_1$, we thus conclude that  for any $t_0>1$,
$$l(t_0)>\lim_{t\to\infty}l(t)= \frac {2}{b_2} \geq l(1)=\frac{2}{\frac{b_2}{2}+\beta+\frac{b_1}{2}},$$
and the proof is therefore complete.\qed

Inspired by \cite{GS}, in the following we reprove Theorem \ref{thm1} by using the  Gagliardo-Nirenberg inequality
(\ref{GNineq}) and some recaling techniques. For the reader's convenience, we restate Theorem \ref{thm1} as the following lemma.

\begin{lem}\label{lemA3}
Let $Q$ be the unique positive radial solution of
(\ref{Kwong}) and suppose  $V_i(x)$ satisfies (\ref{1.12}) for $i=1,\,2$. Then,
\begin{enumerate}
\item [(i)] If $0<b_1<a^*$, $0<b_2<a^*$ and  $0<\beta<\sqrt {(a^*-b_1)(a^*-b_2)}$, then there exists at least one minimizer for (\ref{eq1.4}).

\item [ (ii)] If either $b_1> a^*$ or $b_2>a^*$ or $\beta> \frac{a^*-b_1}{2}+\frac{a^*-b_2}{2}$, then there is no minimizer for (\ref{eq1.4}).
\end{enumerate}
\end{lem}

\noindent\textbf{Proof.} $(i):$ We first note from  the
Gagliardo-Nirenberg inequality (\ref{GNineq}) that
for any $(u_1,u_2)\in \mathcal{X}$ satisfying $\|u_1\|_2^2=\|u_2\|_2^2=1$,
\begin{eqnarray}\label{2:min}
  E_{b_1,b_2,\beta}(u_1,u_2)\geq \sum_{i=1}^2\int_{\R ^2}\Big[\Big(\frac{a^*-b_i}{2}\Big)|
  u_i|^4+V_i(x)|u_i|^2\Big]dx-\beta\int_{\R ^2} |u_1|^2|u_2|^2dx.
\end{eqnarray}
Since $V_i(x)\geq0$ and $0<\beta<\sqrt{(a^*-b_1)(a^*-b_2)}$, one can deduce that $\hat e(b_1,b_2,\beta)\geq 0$
for all $0< b_i< a^*:=\|Q\|^2_2$ ($i=1,\,2$).
Let
$\{(u_{1n},u_{2n})\}\subset \mathcal{M}$ be a minimizing sequence
of problem (\ref{eq1.4}) satisfying
$$\|u_{1n}\|_2^2=\|u_{2n}\|_2^2=1\quad \text{and}\quad
\lim_{n\to\infty}E_{b_1,b_2,\beta}(u_{1n},u_{2n})=\hat e(b_1,b_2,\beta).$$
Taking   $\delta\in(\frac{\beta}{a^*-b_2},\frac{a^*-b_1}{\beta})$,  it then follows from (\ref{2:min}) and   Young's inequality  that
\begin{equation*}
\begin{split}
\hat e(b_1,b_2,\beta)+1&\geq\sum_{i=1}^2\inte \Big(\frac{a^*-b
_i}{2}|u_{in}|^4+V_i(x)|u_{in}|^2\Big)dx-\beta\inte|u_{1n}|^2|u_{2n}|^2dx\\
&\geq\Big(\frac{a^*-b
_1}{2}-\frac{\beta\delta}{2}\Big)\inte|u_{1n}|^4dx+\Big(\frac{a^*-b
_2}{2}-\frac{\beta}{2\delta}\Big)\inte|u_{2n}|^4dx.
\end{split}
\end{equation*}
This implies that  $\{(u_{1n},u_{2n})\}$ is  bounded in $L^4(\R^2)\times L^4(\R^2)$ uniformly w.r.t. $n$, and  it is then easy to deduce that  $\{(u_{1n},u_{2n})\}$ is  bounded uniformly in $\mathcal{X}$. Therefore, by the compactness of Lemma
\ref{2:lem1}, there exist a subsequence of $\{(u_{1n},u_{2n})\}$  and $(u_1,u_2)\in
\mathcal{X}$ such that
\begin{equation*}\begin{split}&(u_{1n},u_{2n})\overset{n}{\rightharpoonup}(u_1,u_2)\quad\text{weakly
in}\ \, \mathcal{X}\,,\\&
(u_{1n},u_{2n})\overset{n}{\to}(u_1,u_2)\quad\text{strongly in}\ \,
L^q(\R^2)\times L^q(\R^2), \end{split}\end{equation*}
where $2\leq q<\infty$. We therefore have $\|u_1\|_2^2=\|u_2\|_2^2=1$ and $
E_{b_1,b_2,\beta}(u_1,u_2)=\hat e(b_1,b_2,\beta)$. This proves the existence of minimizers for the case where $0<b_1<a^*$, $0<b_2<a^*$ and  $0<\beta<\sqrt {(a^*-b_1)(a^*-b_2)}$.

 $(ii):$ Let $\varphi(x)\in C^\infty_0(\R^2)$ be a nonnegative smooth cutoff function such that $\varphi(x)=1$ if $|x|\leq1$ and $\varphi(x)=0$ if $|x|\geq2$. For any $\bar x_0\in\R^2$,  $\tau>0$ and $R>0$, set
 \begin{equation}\label{2:trial}
\phi(x)=A_{R\tau}\frac{\tau}{\|Q\|_2}\varphi\big(\frac{x-\bar
x_0}{R}\big)Q\big(\tau|x-\bar x_0|\big),
 \end{equation}
where $A_{R\tau}>0$  is chosen such that $\int_{\R^2}\phi^2dx=1$. By
scaling, $A_{R\tau}$ depends only on the product $R\tau$. In fact,
using the exponential decay of $Q$ in (\ref{4:exp}), we have
\begin{equation}\label{2:limt1}
\frac{1}{A_{R\tau}^2}=\frac{1}{\|Q\|^2_2}\int_{\R^2}\varphi^2\big(\frac{x}{R\tau}\big)
Q^2(x)dx=1+O((R\tau)^{-\infty})
\quad\text{as}\quad R\tau\to\infty.
\end{equation}
Here we use the notation $f(t) = O(t^{-\infty})$ for a function $f$ satisfying $\lim_{t\to \infty} |f(t)| t^s = 0$ for all $s>0$.
By the exponential decay of $Q(x)$ and the equality (\ref{1:id}), we
have
\begin{equation}\label{2:limt2}
\begin{split}
 \int_{\R^2}|\nabla \phi|^2-\frac{b_i}{2}\int_{\R^2} \phi^4dx
&=\frac{\tau^2}{\|Q\|^2_2}\int_{\R^2}|\nabla
Q|^2-\frac{b_i\tau^2}{2\|Q\|^4_2}\int_{\R^2} Q^4dx+O((R\tau)^{-\infty})\\
&=\Big(1-\frac{b_i}{\|Q\|_2^2}\Big)\tau^2+O((R\tau)^{-\infty})
\quad\text{as}\quad R\tau\to\infty.
\end{split}\end{equation}
On the other hand, since the function $x\mapsto
V_i(x)\varphi^2(\frac{x-\bar x_0}{R})$ is bounded and has compact
support,  we obtain  that
\begin{equation}\label{2:limt6}
\lim_{\tau\to\infty}\int_{\R^2}V_i(x)\phi^2(x)dx=V_i(\bar x_0)
\end{equation}
holds for almost every $\bar x_0\in\R^2$, where $i=1,\, 2$.

Suppose that  $b_1>a^*$. Choosing
$\eta(x)\in C_0^\infty(\R^2)$ such that $\int_{\R^2}\eta^2(x)dx=1$,
we then derive that
\begin{equation*}
\int_{\R^2}\phi^2\eta^2dx\leq \sup_{x\in
R^2}\eta ^2(x)\int_{\R^2}\phi^2dx=\sup_{x\in R^2}\eta ^2(x)<\infty,
\end{equation*}
where $\phi>0$ is as in (\ref{2:trial}). Together with (\ref{2:limt2}) and  (\ref{2:limt6}), this estimate  then yields that
\begin{equation*}\begin{split}E_{b_1,b_2,\beta}(\phi,\eta)
\leq  \inte|\nabla\phi|^2-\frac{b_1}{2}\int_{\R^2} \phi^4dx +C\leq \frac{a^*-b_1}{a^*}\tau^2+C,
\end{split}\end{equation*}
which implies that
$$\hat e(b_1,b_2,\beta)\leq \lim_{\tau\to \infty}E_{b_1,b_2,\beta}(\phi,\eta)=-\infty.$$
Similarly, this estimate is still true  if $b_2>a^*$. Therefore,  $\hat e(b_1,b_2,\beta)$ does not admit any minimizer.

Finally, if  $\beta> \frac{a^*-b_1}{2}+\frac{a^*-b_2}{2}$, we also deduce  from (\ref{2:limt2}) and (\ref{2:limt6}) that
\begin{equation*}
\hat e(b_1,b_2,\beta)\leq \lim_{\tau\to \infty}E_{b_1,b_2,\beta}(\phi,\phi)=-\infty,
\end{equation*}
which also implies the
non-existence of minimizers. This completes the proof of Lemma \ref{lemA3}.\qed\\

\noindent {\bf Acknowledgements:} This work was  supported
    by National Natural Science Foundation of China (11322104, 11271360, 11471331, 11501555) and National Center for Mathematics and Interdisciplinary Sciences.\\

\noindent {\bf Conflict of interest statement:} The authors declare that they have no conflict of interest.

\end{document}